\documentclass[12pt,letterpaper]{article}

\usepackage{amsmath}
\usepackage{amssymb}
\usepackage{amsthm}
\usepackage{enumitem}
\usepackage{graphicx}
\usepackage{hyperref}

\newcommand\bs{\backslash}
\newcommand{\C}{\mathbb{C}}

\renewcommand{\O}{\mathcal{O}}
\newcommand{\oil}{of index~$\ell$}
\newcommand{\p}{\mathfrak{p}}
\newcommand{\Q}{\mathbb{Q}}
\newcommand{\R}{\mathbb{R}}
\newcommand{\Z}{\mathbb{Z}}
\newcommand{\GL}{\mathrm{GL}}
\newcommand{\SL}{\mathrm{SL}}

\newcommand{\vor}{Vo\-ro\-noi}

\newtheorem{corollary}{Corollary}
\newtheorem{lemma}{Lemma}
\newtheorem{proposition}{Proposition}
\newtheorem{theorem}{Theorem}

\theoremstyle{definition}
\newtheorem{definition}{Definition}
\newtheorem{algorithm}{Algorithm}

\theoremstyle{remark}
\newtheorem*{remark}{Remark}

\newtheorem*{examples}{Examples}

\begin{document}
\date{March 21, 2024}
\title{Tempered Perfect Lattices in the Binary Case}
\author{Erik Bahnson \and Mark McConnell \and Kyrie McIntosh
  \thanks{Department of Mathematics, Princeton University.
    \\ E-mails: \texttt{ebahnson@alumni.princeton.edu},
    \texttt{markwm@princeton.edu},
    \texttt{kyriem@alumni.princeton.edu}.  }}

\maketitle

\begin{abstract}
  A new algorithm for computing Hecke operators for $\SL_n$ was
  introduced in \cite{MM20}.  The algorithm uses \emph{tempered
  perfect lattices}, which are certain pairs of lattices together with
  a quadratic form.  These generalize the perfect lattices of
  \vor\ \cite{Vor}.  The present paper is the first step in
  characterizing tempered perfect lattices.  We obtain a complete
  classification in the plane, where the Hecke operators are for
  $\SL_2(\Z)$ and its arithmetic subgroups.  The results depend on the
  class field theory of orders in imaginary quadratic number fields.
\end{abstract}

%%%%%%%%%%%%%%%%%%%%%%%%%%%%%%%%%%%%%%%%

\section{Introduction}

Lagrange and Gauss knew that a plane lattice~$L$ with a quadratic form
could have interesting number-theoretic structure.  This paper
describes the number-theoretic structure of pairs of lattices $M
\subset L$.  Our methods are elementary, but class field theory is
needed to tell the whole story.

The motivation is to study a new algorithm \cite{MM20} for computing
Hecke operators.  Let $\mathbb{G}$ be the algebraic group~$\SL_n$ of a
finite-dimensional division algebra (possibly commutative) over~$\Q$.
Let~$G$ be the real points of~$\mathbb{G}$, and let~$X$ be the
symmetric space for~$G$.  Let $\Gamma\subset G$ be an arithmetic
subgroup.  We study the cohomology $H^i(\Gamma\bs X; \rho)$ with
coefficients in suitable finite-dimensional coefficient
systems~$\rho$.  A family of Hecke operators acts on this cohomology.
Many authors have given algorithms for computing these Hecke operators
in different settings and for different ranges of the degree~$i$.  (A
short list relevant to this paper is \cite{Manin, Stein, AR, Gun}.)
The new algorithm in~\cite{MM20} works in complete generality in the
setting we have described---$\SL_n$ of any division algebra, with
arbitrary coefficients~$\rho$, and in all degrees~$i$.  See
\cite[\S2]{GM21} for an exposition of \cite{MM20} in the most common
case $G = \SL_n(\R)$.

The new algorithm uses \emph{tempered perfect lattices}.  These are
certain pairs of lattices $M\subset L$ together with a quadratic form.
They generalize the perfect forms introduced by \vor~\cite{Vor}.
Section~\ref{sec:bg} will define these terms.

The definition of a tempered perfect lattice uses only Euclidean
geometry.  When we run the algorithm of~\cite{MM20}, however, the data
shows that tempered perfect lattices also depend on number theory.
This paper is the first step in understanding and characterizing
tempered perfect lattices.  For this first step, we restrict to
$G=\SL_n(\R)$.  In fact, we restrict to $\SL_2(\R)$, because the
classical theory of binary quadratic forms is so rich.

For $\SL_2(\Z)$, there is essentially one Hecke operator $T_\ell$ at
each prime~$\ell$, and these~$T_\ell$ generate the algebra of all
Hecke operators.  In any approach to Hecke operators for $\SL_2(\Z)$,
$T_\ell$ is defined using lattices~$L$ and sublattices $M\subset L$ of
index~$\ell$.  We assume $M\subset L$ has prime index~$\ell$
throughout the paper.

We now briefly describe the different families of tempered perfect
lattices for $n=2$.  We refer to figures that come later in the paper.
In the figures, points of~$L$ are dots~$\cdot$, points of~$M$ are
circled dots~$\odot$, and the quadratic form is the standard dot
product.  The \emph{temperament}~$\tau$ is defined to be the ratio of
the lengths of the shortest non-zero vectors in~$M$ and~$L$,
respectively.  This is the ratio of the radii of the two circles in
each figure.  When there are $s$ pairs $\pm v$ of vectors of $L-M$ on
the inner circle, and $s'$ pairs of vectors of~$M$ on the outer
circle, we call the form an \emph{$s$-and-$s'$ form}.

The only perfect plane lattice, up to rotation and scaling, is the
hexagonal lattice.  One family of tempered perfect lattices is the
\emph{$3$-and-$3$ forms}, where~$L$ and~$M$ are both hexagonal
lattices, a snowflake within a snowflake.  Figure~\ref{fig:exEis7}
shows the 3-and-3 form with index $\ell=7$.  We identify~$L$ with the
Eisenstein integers, the ring $\Z[\omega]$ where $\omega = e^{2\pi i
  /3}$.  Figure~\ref{fig:exEis7} arises because the prime~7 splits as
$(2-\omega)(2-\bar{\omega})$.  Here $2-\omega$ is one of the three
pairs of shortest vectors in~$M$, and all these vectors have length
$\tau = \sqrt{7}$.
%% The index is~7 because $M$ comes from~$L$ by scaling by~$\sqrt{7}$ in
%% both the~$x$ and~$y$ directions, followed by a rotation.

Figure~\ref{fig:exEis11} is a \emph{$3$-and-$1$ form}.  Here~$L$ is
again $\Z[\omega]$.  The shortest non-zero vector $v\in M$ still has
length~$\sqrt{7}$, but in any lattice basis $\{v,w\}$ of~$M$, $w$ is
strictly longer than~$\sqrt{7}$.  Thus~$M$ has less symmetry than~$L$.
The 3-and-1 forms are determined partly by number theory and partly by
geometry.  Theorem~\ref{thm:eis} and Algorithm~\ref{alg:eis} give
their precise characterization.

A lattice is \emph{well rounded} \cite{Ash84} if the non-zero lattice
vectors of minimal length contain a basis.  In any well-rounded plane
lattice (besides the hexagonal lattice), the shortest non-zero vectors
comprise two pairs $\pm v, \pm w$, and $\{v,w\}$ is a lattice basis.
When~$L$ and~$M$ are both well rounded with two pairs of minimal
vectors each, then $M\subset L$ is a \emph{$2$-and-$2$ form}.

A familiar well-rounded lattice is the square lattice $L$, which may
be identified with the Gaussian integers $\Z[i]$.  In
Figure~\ref{fig:exGauss17}, $M$ is $1+4i$ times~$L$.  Then~$M$ is also
a square lattice, sitting in~$L$ at an angle.  This is a 2-and-2 form.
The index is the norm $|1+4i|^2 = 17$.  The temperament is $\tau =
\sqrt{17}$.
%% since both of the edges in the square grid in~$M$ have length $|z|
%% = \sqrt{17}$.

The 2-and-2 forms become more interesting when~$M$ has a different
shape from~$L$.  Figure~\ref{fig:ex1155} shows the first example,
which occurs for $\ell=23$.  The angle between the vectors in a
well-rounded basis of~$L$ (the points on the inner circle) is not
equal to the angle for~$M$ (the outer circle).  For the first time,
$\tau^2$ is not an integer: it is $\frac{17\cdot23}{19}$.  We soon
discover that class field theory determines 2-and-2 forms.  $L$
and~$M$ are fractional ideals in the quadratic number field~$K =
\Q(\sqrt{-1155})$, where $1155 = 3\cdot 5 \cdot 7 \cdot 11$.  The
class group of~$K$ has order eight, and the seven primes
$3,5,7,11,17,19,23$ encountered in this example lie in exactly the
seven non-trivial classes of the class group (Table~\ref{tab:1155}).
Well-rounded fractional ideals of~$K$ lie in the classes for~$17$
and~$19$ and in no others.  The class of~$23$ times the class of~$19$
is the class of~$17$, and this is why~$\tau^2 = \frac{17\cdot23}{19}$
has a 2-and-2 form.

We now summarize our main theorems characterizing 2-and-2 forms,
theorems that are stated and proved in Sections~\ref{sec:22asFuncEll}
and~\ref{sec:22asFuncD}.  We show that~$L$ and~$M$ are proper ideals
in an order~$\O$ of discriminant $D<0$ in a quadratic number field.
$\O$ depends on~$\ell$, $L$, and~$M$, but it is the same for~$M$ as
for~$L$.  Hence the 2-and-2 forms are controlled by arithmetic in the
class group of~$\O$.  Theorem~\ref{thm:finManyDForEll} shows that only
finitely many orders~$\O$ will yield 2-and-2 forms for a given
index~$\ell$, and it gives a bound on the discriminant~$D$ in terms
of~$\ell$.  Theorem~\ref{thm:infManyEll} is a converse: if an
order~$\O$ has a 2-and-2 form for a certain~$\ell$, the same~$\O$ will
yield 2-and-2 forms for infinitely many~$\ell'$, those coming from the
Chebotaryov Density Theorem applied to the class of~$\ell$ in~$\O$.  In
Section~\ref{subsec:ruleOfThree}, we show why the discriminants~$D$
for 2-and-2 forms will always be a product of comparatively small
prime factors.

There are many directions for future work growing out of this paper.
One is to classify tempered perfect lattices for~$\Gamma \subseteq
\SL_2(\O_F)$, where~$\O_F$ is the ring of integers of a number
field~$F$.  This would connect with the literature on Hecke operators
for Bianchi modular forms and Hilbert modular forms.  (As a short
list, see \cite{AGMY, GunMY} and the references therein.)  A second
direction is to classify tempered perfect forms for~$\Gamma \subseteq
\SL_n(\Z)$ for $n>2$, especially for $n=3$ and~$4$.  This would
involve pairs of projective configurations, as in \cite{MM89}.

%%%%%%%%%%%%%%%%%%%%%%%%%%%%%%%%%%%%%%%%

\section{Background} \label{sec:bg}

\subsection{Lattices and pairs of lattices}

We view $\R^2$ as a space of row vectors.  When $\{v,w\}$ is a basis
of~$\R^2$, a \emph{lattice}~$L$ is the set of $\Z$-linear combinations
of~$v$ and~$w$.  We call $\{v,w\}$ a \emph{$\Z$-basis} of~$L$.
%% A lattice is a discrete subgroup of~$\R^2$ and is isomorphic to~$\Z^2$.

Let~$L_0$ be the standard copy of~$\Z^2$ with basis $\{(1,0),
(0,1)\}$.  Let $G = \GL_2(\R)$ (not~$\SL_2(\R)$ as in the
Introduction).  The subgroup of~$G$ that stabilizes~$L_0$ is $\Gamma_0
= \GL_2(\Z)$.  Every lattice~$L$ has the form $L_0\,g$ for $g\in G$.
%% Lattices $L_0\,g$ and $L_0\,g'$ are equal if and only if $g' = \gamma g$ for $\gamma\in\Gamma_0$.
When we have a specific~$g$ in mind, we say the rows of~$g$ are a
\emph{distinguished basis} of~$L$.
%% The quotient $Y = \Gamma_0\bs G$ is the \emph{space of lattices} in $\R^2$.

A \emph{homothety} is multiplication of~$\R^2$ by a positive real
scalar.  The homotheties form a group $\R_+$.

By $(L,M)$ we denote a pair $M \subset L$ where~$L$ is a lattice
and~$M$ is a sublattice of finite index.  $M$ has the form~$Lg$ for
some integer matrix~$g$ with $[L:M] = |\det g|$.  In the special case
$M_0 = L_0\,g$ for an integer matrix~$g$, the group $\Gamma = \Gamma_0
\cap g^{-1} \Gamma_0\,g$ is the common stabilizer in~$G$ of~$L_0$
and~$M_0$, the \emph{stabilizer of the pair} $(L_0,M_0)$.
This~$\Gamma$ has finite index in~$\Gamma_0$.

\subsection{Pairs of lattices of prime index} \label{subsec:indexEll}

Throughout the paper, let~$\ell$ be a prime number.  A 
lattice~$L$ contains exactly $\ell+1$ sublattices~$M$ of index~$\ell$,
namely
\begin{equation} \label{heckeMatsEll}
\begin{bmatrix} \ell & 0 \\ 0 & 1 \end{bmatrix} L, \qquad
\begin{bmatrix} 1 & k \\ 0 & \ell \end{bmatrix} L
\end{equation}
where~$k=0,\dots,\ell-1$.
%% The proof is standard \cite[VII.5.2]{Se} \cite[Prop.~7.2]{Bu}.
It is easy to show that~$\Gamma_0$ acts transitively on the
sublattices of~$L_0$ of index~$\ell$.

Suppose $M_0$ has basis $\{(1,0), (0,\ell)\}$.  The stabilizer
of~$(L_0,M_0)$ is the classical congruence subgroup
\[
\Gamma_0(\ell) = \left\{ \left. \begin{bmatrix} a & b \\ c &
  d \end{bmatrix} \in \Gamma_0 \, \right| \, b\equiv 0 \bmod \ell
\right\}
\]
whose index in~$\Gamma_0$ is $\ell+1$.

%% The quotient $Y(\ell) = \Gamma_0(\ell)\bs G$, or with $\Gamma_0(\ell)$
%% replaced by a conjugate, is the space of index-$\ell$ pairs of
%% lattices in $\R^2$.  This is the space of lattices with
%% \emph{level~$\ell$ structure}.

An element $v\in L$ is \emph{primitive} if~$v$ generates $(\R v)\cap L$.

\begin{lemma} \label{uniqueSublatt}
  A given primitive vector in a lattice~$L$ is contained in one and
  only one sublattice~$M$ of index~$\ell$.
\end{lemma}

\begin{proof}
  It suffices to consider a primitive $(x,y)\in L_0$.  Take integers
  $s,t$ with $sx+ty=1$.  Multiplying by $\begin{bmatrix} x & y \\ -t &
    s \end{bmatrix}^{-1} \in \Gamma_0$, we may assume the primitive
  vector is $(1,0)$.  In~\eqref{heckeMatsEll}, only $\begin{bmatrix} 1
    & 0 \\ 0 & \ell \end{bmatrix}$ has $(1,0)$ in the $\Z$-span of its
  rows.
\end{proof}
  
\subsection{Quadratic forms} \label{subsec:quforms}

Throughout the paper, let $\langle (x_1,y_1), (x_2,y_2) \rangle
= \begin{bmatrix} x_1 & y_1 \end{bmatrix} \begin{bmatrix}
  1&0\\0&1 \end{bmatrix} \begin{bmatrix} x_2 \\ y_2 \end{bmatrix}$ be
the standard dot product on~$\R^2$.
%%% We say $\begin{bmatrix} 1&0\\0&1 \end{bmatrix}$ is the matrix of
%%% the inner product.
The subgroup of~$G$ preserving~$\langle\,,\,\rangle$ is the group
$\mathrm{O}_2(\R)$ of orthogonal matrices.  Let $X =
G/\mathrm{O}_2(\R)\R_+$; this is the symmetric space from the
Introduction.

Each lattice $L = L_0\,g$ carries the inner product
$\langle\,,\,\rangle$.  We can view this from a second perspective:
not a fixed inner product on a space of lattices, but a fixed lattice
carrying a space of inner products.  If we change coordinates
by~$g^{-1}$, then~$L$ becomes the fixed lattice~$L_0$, and the matrix
of the inner product becomes $A = g\begin{bmatrix}
1&0\\0&1 \end{bmatrix} g^T = g g^T$, where~${}^T$ is the transpose.
This defines a quadratic form $\begin{bmatrix} x & y \end{bmatrix}
A \begin{bmatrix} x \\ y \end{bmatrix}$ on~$L_0$.
%% By abuse of notation, we will also say~$A$ is a quadratic form.
In this paper, quadratic forms will always be positive definite.
%%% $A$ is a \emph{binary} quadratic form because there are only two
%%% variables $(x,y)$.

We consider two quadratic forms to be equivalent if they differ by a
homothety.  Our~$X$ is the space of quadratic forms modulo
homotheties.  Going back to the lattice perspective, $X$ is the space
of lattices with a distinguished basis, in which two lattices are
considered equivalent if one can be carried to the other by an
orthogonal transformation and a homothety.

\subsection{Reduced bases} \label{subsec:redBases}

The \emph{minimal vectors} of~$L$ are the $v\in L - \{0\}$ where
$\|v\|^2 = \langle v,v\rangle$ attains its minimal value.  Minimal
vectors occur in pairs $\pm v$.

The lattice~$L$ has a \emph{Minkowski reduced basis}, or just
\emph{reduced basis}, constructed as follows \cite[(2.7.7)]{BoSh}.
Let~$v \in L$ be a minimal vector.  Let $w\in L$ have minimal length
among vectors for which $\{v,w\}$ is linearly independent over~$\R$.
Then $\{v,w\}$ is a $\Z$-basis of~$L$, and the angle~$\theta$
between~$v$ and~$w$ satisfies $\pi/3 \leqslant \theta \leqslant
2\pi/3$.
%% \cite[pp.~77--82]{Se}
(We emphasize that this is for $n=2$.  For $n \geqslant 5$, naively
picking the shortest lattice vectors in successive dimensions
$1,2,\dots,n$ may not produce a lattice basis \cite[(15.10.1)]{CS}.)

\begin{figure}
  \begin{center}
    \includegraphics[scale=0.1]{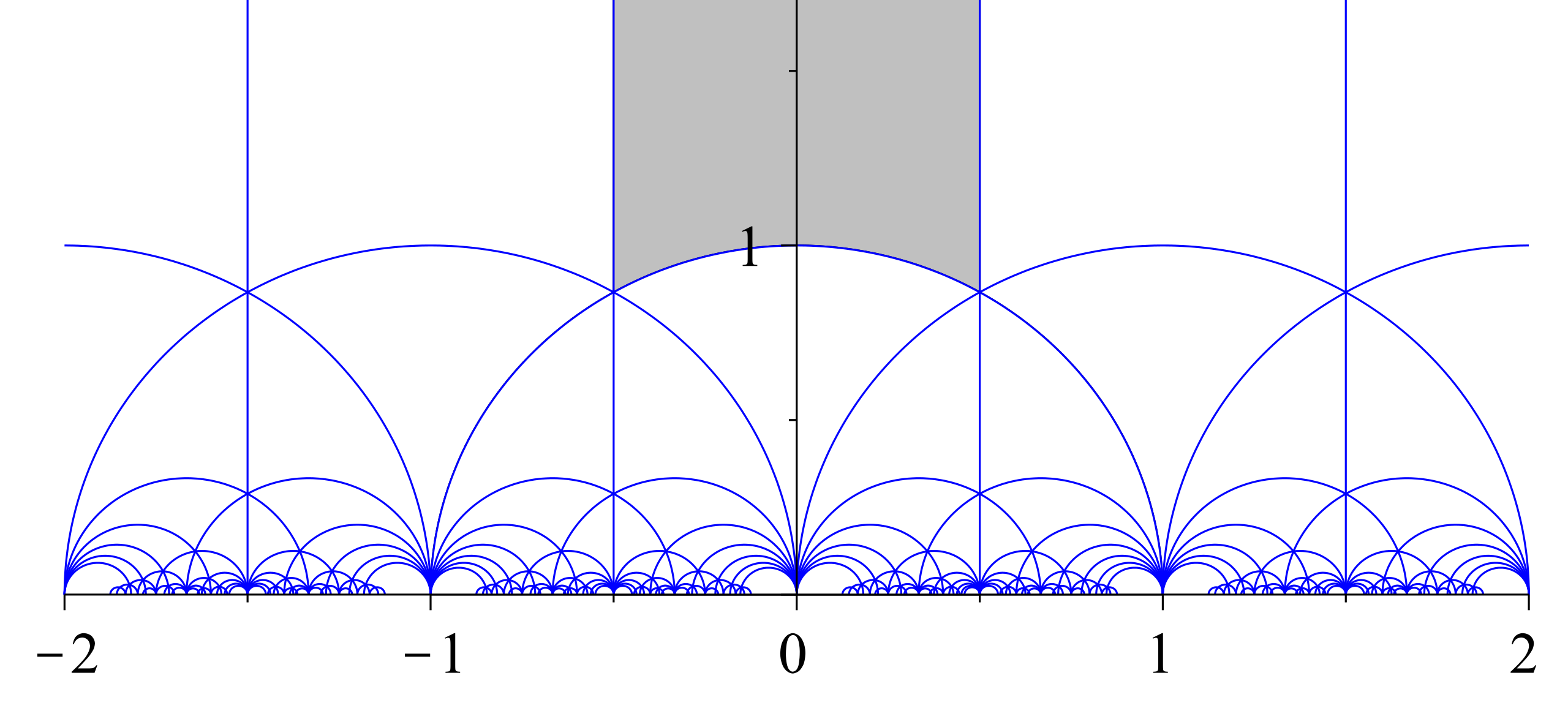}
  \end{center}
  \caption{A reduced basis $\{v,w\}$ has $v=(1,0)$ and $w$ in
    the shaded region.}
  \label{fig:sl2zFundDom}
\end{figure}

In Figure~\ref{fig:sl2zFundDom}, we reproduce the standard picture.
Rescale~$L$ by a homothety so that the minimal vector~$v$ has
length~1, and rotate to carry~$v$ to $(1,0)$.  Then the lattice basis
$\{v,w\}$ is reduced if and only if~$\pm w$ is in the shaded region.
This region is a fundamental domain for $\SL_2(\Z)\bs X$.

When $n=2$, a lattice is well rounded if, among its minimal vectors,
there are two vectors forming a $\Z$-basis of the lattice.  This
$\Z$-basis is Minkowski reduced, and its~$w$ is on the bottom circular
edge of the shaded region in Figure~\ref{fig:sl2zFundDom}.  When the
angle is exactly $\pi/3$ or $2\pi/3$, we get the hexagonal lattice.
No lattice for $n=2$ has four or more pairs of minimal vectors.

\subsection{Tempered perfect forms}

We refer to \cite{MM20} for rather lengthy technical material on
temperaments, the well-tempered complex and its cell decomposition,
and critical temperaments, as well for how the complex is used to
compute Hecke operators.  This relies on \cite{Ash84} for sets of
weights, the well-rounded retract, and perfect forms.  Assuming this
technical background, a tempered perfect lattice is by definition a
$0$-dimensional cell in the well-tempered complex.  A tempered perfect
form is a quadratic form on~$L_0$ corresponding to a tempered perfect
lattice.

We now give an equivalent, self-con\-tained definition of a tempered
perfect form.

\begin{definition} \label{def:tauperf}
  Let~$\ell$ be a prime number.  Let $L_0 = \Z^2$ and let $M_0$ be a
  sublattice of~$L_0$ of index~$\ell$.  A positive definite quadratic
  form~$A$ on~$L_0$ is a \emph{tempered perfect form} if there exist
  non-empty subsets $S \subset L_0 - M_0$ and $S' \subset M_0 - \{0\}$
  such that

  (i) $v A v^T = 1$ for all $v\in S$;

  (ii) $v A v^T > 1$ for all $v\in L_0 - M_0$ with $v\notin S$;

  (iii) $v A v^T$ have a common value~$u$ for all $v\in S'$, with $u
  \geqslant 1$;

  (iv) $v A v^T > u$ for all $v\in M_0 - \{0\}$ with $v\notin
  S'$;

  (v) if~$B$ is any quadratic form satisfying (i)--(iv) with common
  value $u'$, then $B=A$ and $u=u'$.
\end{definition}

The \emph{temperament} of a tempered perfect form is $\tau = \sqrt{u}$.

We use the constant~$1$ in (i)--(ii) to normalize the form in its
homothety class.  If $\tau=1$, we recover \vor's definition of perfect
form \cite[\S\S1--5]{Vor}, in which $S\cup S'$ is the (single) set of
minimal vectors.

Similarly, $L$ is a \emph{tempered perfect lattice} if there is a
sublattice $M \subset L$ \oil, and non-empty sets $S \in L - M$, $S'
\subset M - \{0\}$, such that (i)--(v) hold taking~$A$ to be the
standard dot product.  We will go freely back and forth between the
perspectives of tempered perfect lattices and forms for the rest of
the paper.
%% While the computation of tempered perfect lattices and forms is the
%% same, the objects have different visual interpretations.

\begin{definition} \label{def:sSp}
  An \emph{$s$-and-$s'$ form} is a tempered perfect form where~$s$
  (resp.,~$s'$) is the number of pairs $\pm v$ in~$S$ (resp.,~$S'$).
\end{definition}

Since well-rounded lattices have either two or three pairs of minimal
vectors,

\begin{definition}
  A \emph{doubly well-rounded form} is an $s$-and-$s'$ form with
  $s\geqslant 2$ and $s' \geqslant 2$.
\end{definition}

%%%%%%%%%%%%%%%%%%%%%%%%%%%%%%%%%%%%%%%%

\section{Structure of tempered perfect forms}

\subsection{Rationality of tempered perfect forms}

Write a tempered perfect form as $ax^2 + bxy + cy^2$ for $(x,y)\in
L_0$.  The matrix of the form is $\begin{bmatrix} a & b/2 \\ b/2 &
  c \end{bmatrix}$.  All we know at first is that $a,b,c\in\R$.  The
next proposition generalizes \vor's result that perfect forms are
rational \cite[\S5]{Vor}.

\begin{proposition} \label{oldconj1}
  For a tempered perfect form $ax^2 + bxy + cy^2$, the quantities~$a$,
  $b$, $c$, and $\tau^2$ are in~$\Q$.
\end{proposition}

\begin{proof}
  We use the notation of Definitions~\ref{def:tauperf}
  and~\ref{def:sSp}.  Let $(x_i,y_i)\in\Z^2$ be the vectors in~$S$
  followed by~$S'$, with only one from each pair~$\pm v$.  Parts~(i)
  and~(iii) of the definition state
  \begin{equation} \label{rationalityProp}
  \begin{bmatrix}
    x_1^2 & x_1 y_1 & y_1^2 & 0 \\
    \vdots & \vdots & \vdots \\
    x_s^2 & x_s y_s & y_s^2 & 0 \\
    x_{s+1}^2 & x_{s+1} y_{s+1} & y_{s+1}^2 & -1 \\
    \vdots & \vdots & \vdots & \vdots \\
    x_{s+s'}^2 & x_{s+s'} y_{s+s'} & y_{s+s'}^2 & -1
  \end{bmatrix}
  \begin{bmatrix}
    a \\ b \\ c \\ u
  \end{bmatrix} =
  \begin{bmatrix}
    1 \\ \vdots \\ 1 \\ 0 \\ \vdots \\ 0
  \end{bmatrix}.
  \end{equation}
  Part~(v) states there is a \emph{unique} solution for $a$, $b$, $c$,
  $u$ when we regard the $x_i$ and $y_i$ as given.  Parts~(ii)
  and~(iv) do not affect the uniqueness statement, since they are
  strict inequalities (open conditions).  Thus the system of linear
  equations~\eqref{rationalityProp} must have a unique solution.  The
  coefficient matrix on the left has entries in~$\Z$, as does the
  vector on the right.  Hence the solution has entries in~$\Q$.  By
  definition, $\tau^2 = u$.
\end{proof}

Since the matrix in~\eqref{rationalityProp} has rank~$4$,

\begin{corollary} \label{sPlusSp}
  In an $s$-and-$s'$ form, $s + s' \geqslant 4$.
\end{corollary}

\subsection{The boundary temperaments} \label{bdryTau}

The temperament $\tau=1$ is special.  The only perfect lattice for
$n=2$ is the hexagonal lattice, so perfect forms for $\tau=1$ have
three pairs of minimal vectors.  Either one or none of these pairs
will be in~$M$.  By Corollary~\ref{sPlusSp}, no form for $\tau=1$ is
tempered perfect.

The temperament $\tau=1$ is a boundary case because $\tau\geqslant 1$
by~(iii).  The next lemma shows that $\tau=\ell$ is the other boundary
case.

\begin{lemma}
  In any tempered perfect form, $\tau \leqslant \ell$.
\end{lemma}

\begin{proof}
  Suppose $\tau>\ell$.  Let~$v$ be the shortest non-zero vector
  in~$L$.  Then $vAv^T = 1$ by~(i), so $(\ell v)A(\ell v)^T = \ell^2 <
  u$.  But $\ell v \in M$, so $(\ell v)A(\ell v)^T = \ell^2 \geqslant
  u$ by (iii)--(iv), a contradiction.
\end{proof}

Consider an $s$-and-$s'$ form in the boundary case $\tau=\ell$.  If
$v\in S$, then $\ell v \in S'$ since $u=\ell^2$.  Hence $s\leqslant
s'$.  We cannot have $s\geqslant 2$, for then $s'\geqslant 2$, and~$M$
would contain a well-rounded sublattice of index~$\ell^2$; the index
$[L:M]$ is prime, hence~$M$ would contain some other vector shorter
than squared-length~$\ell^2$, contradicting the definition of~$u$.

We can have tempered perfect forms with $\tau=\ell$ with $s=1$.
Corollary~\ref{sPlusSp} implies these are 1-and-3 forms, where~$M$ is
the hexagonal lattice.

Turning away from the boundary cases, we now suppose $1<\tau<\ell$.
All the minimal vectors for~$L$ are in~$S$ but not in~$S'$.  The
shortest non-zero vector in~$M$ is always in~$S'$.  Thus an
$s$-and-$s'$ form has $s\geqslant 1$ and $s'\geqslant 1$.  For doubly
well-rounded forms, $S$ contains a Minkowski reduced basis of~$L$, and
$S'$ contains a Minkowski reduced basis of~$M$.

\subsection{Enough minimal vectors implies tempered perfect}

We will need the following converse to Proposition~\ref{oldconj1}.

\begin{proposition} \label{sPlusSpConverse}
If a rational quadratic form $A$ satisfies parts~(i)-(iv) of
Definition~\ref{def:tauperf} with $s + s' \geqslant 4$, then it also
satisfies part~(v) and is thus a tempered perfect form.
\end{proposition}
\begin{proof}
It suffices to prove that~\eqref{rationalityProp} has a unique
solution for $a,b,c,u$ when $s + s' \geqslant 4$.

From Section~\ref{subsec:redBases}, $s\leqslant 3$ and $s'\leqslant
3$.  First suppose $s = 3$.  Then~$L$ is the hexagonal lattice.  After
a change of coordinates, $A = x^2 - xy + y^2$.  This determines~$a$,
$b$, and~$c$, and~$u$ is determined from~(iii) by any one of the
vectors in~$S'$.  When $s' = 3$, the proof is similar, interchanging
the roles of~$S$ and~$S'$.

The remaining case is $s, s'\geqslant 2$.  From Section~\ref{bdryTau}, $1 <
u < \ell^2$.  We may change coordinates so that $L=L_0$ and~$S = \{\pm
v_1, \pm v_2\}$ for $v_1 = (1,0)$ and $v_2 = (0,1)$.  Let $v_3 = (x_3,
y_3)$, $v_4 = (x_4, y_4) \in S'$ with $v_4 \ne \pm v_3$.  Because
$u<\ell^2$, $\ell v_1$ and $\ell v_2$ are not in~$M$, so $x_3$,
$y_3$, $x_4$, $y_4$ are all non-zero.  Let $A = \begin{bmatrix}a & b/2 \\ b/2 &
  c \end{bmatrix}$.  We use $v_1$ and $v_2$ to deduce $a = c = 1$.
Then $v_3 A v_3^T = v_4 A v_4^T = u$.  Since $\{v_1,v_2\}$ is
Minkowski reduced, $|b|\leqslant 1$.  As in~(v), suppose $B
= \begin{bmatrix} a' & b'/2 \\ b'/2 & c' \end{bmatrix}$ is an
arbitrary quadratic form satisfying (i)--(iv).  Then $a'= a = 1$, $c'
= c = 1$, $|b'|\leqslant 1$, and $v_3 B v_3^T = v_4 B v_4^T = u'$.  We
have
\begin{align}
  u  &= x_3^2 + y_3^2 + x_3y_3b   \notag \\
  u  &= x_4^2 + y_4^2 + x_4y_4b   \label{fourUUp} \\
  u' &= x_3^2 + y_3^2 + x_3y_3b'  \notag \\
  u' &= x_4^2 + y_4^2 + x_4y_4b'. \notag
\end{align}
Subtracting,
\[
u-u' = x_3y_3(b-b') = x_4y_4(b-b').
\]
If $b=b'$, then $u=u'$ and we are done.  So suppose $b\ne b'$.  Then
$x_3y_3 = x_4y_4$.  From the first two equations in~\eqref{fourUUp},
$x_3^2+y_3^2 = x_4^2+y_4^2$.  Changing the sign of~$b$ if necessary,
$x_3>0$ and $y_3>0$.  Considering the way a circle at the origin meets
the hyperbola $xy=C$ for $C>0$, we see there are only two
possibilities, $(x_3,y_3)=(x_4,y_4)$ and $(x_3,y_3)=(y_4,x_4)$.  We
have excluded the first possibility, so $x_3 = y_4 > 0$ and $x_4 = y_3
> 0$.  Interchanging~$v_3$ and~$v_4$ if necessary, we may assume $x_3
> y_3 > 0$.

Since $\{v_3,v_4\}$ is a Minkowski reduced basis of~$M$, we have $\det
\begin{bmatrix} x_3 & y_3 \\ y_3 & x_3 \end{bmatrix} =
(x_3+y_3)(x_3-y_3) = \ell$.  Since $\ell$ is prime, $x_3+y_3 = \ell$
and $x_3-y_3 = 1$.  Thus~$\ell$ must be odd, $x_3 = (\ell+1)/2$, and
$y_3 = (\ell-1)/2$.  But then $v_3 - v_4 = (1,-1) \in M$.  Hence
\[
\begin{bmatrix} 1 & -1 \end{bmatrix}
\begin{bmatrix} 1 & b/2 \\ b/2 & 1 \end{bmatrix}
\begin{bmatrix} 1 \\ -1 \end{bmatrix}
> u.
\]
We deduce $2-b > u = x_3^2+y_3^2+x_3y_3b$, or $2 > \ell^2/2 + 1/2 +
(\ell^2+3)b/4$, or $b < \frac{6-2\ell^2}{3+\ell^2} =
\frac{12}{3+\ell^2}-2$.  Because $|b|\leqslant 1$, we have
$\frac{12}{3+\ell^2}-2 \geqslant -1$, which forces $\ell=3$, $b=-1$,
$v_3=(2,1)$, and $v_4=(1,2)$.
%% We did not change the sign of~$b'$ above, and one finds $b'=1$.
Changing coordinates so that the quadratic form becomes the standard
dot product, we have arrived at the 3-and-3 form where~$L$ and~$M$ are
both hexagonal lattices and $[L:M]=3$, which is tempered perfect.
\end{proof}

%%%%%%%%%%%%%%%%%%%%%%%%%%%%%%%%%%%%%%%%

\section{Background on binary quadratic forms} \label{sec:eulerLagrGauss}

Proposition~\ref{oldconj1} shows that, up to homothety, a tempered
perfect form is a binary quadratic form $ax^2 + bxy + cy^2$ with
$a,b,c\in\Z$ and $\gcd(a,b,c) = 1$.  The form can also be written
$(a,b,c)$.  A form with $\gcd(a,b,c) = 1$ is \emph{primitive}.  The
\emph{discriminant} of $ax^2 + bxy + cy^2$ is $D = b^2 - 4ac$.  As we
have said, all forms in this paper are positive definite; this implies
$a > 0$, $c > 0$, and $D < 0$.
%% In this section, we briefly review binary quadratic forms.  Modern
%% sources for this material are \cite{BoSh, Cohn, Bu, Cox}.

\subsection{Binary quadratic forms and lattices}

As in Section~\ref{subsec:quforms}, $g \mapsto gg^T$ maps a lattice
with distinguished basis to a quadratic form.  The matrix of $ax^2 +
bxy + cy^2$ is $\begin{bmatrix} a & b/2 \\ b/2 & c \end{bmatrix}$.
Dividing by~$a$ and solving $\begin{bmatrix} 1 & b/2a \\ b/2a &
  c/a \end{bmatrix} = g g^T$ for~$g$, we find one solution
\begin{equation} \label{choleskygamma}
  g = \begin{bmatrix} 1 &
    0 \\ b/2a & \sqrt{-D}/2a \end{bmatrix}.
\end{equation}
(The general solution is $gk$ for~$k$ orthogonal.)

A special property of $n=2$ is that these lattice bases can be
understood using complex numbers.  Identify $\R^2$ with $\C$ by $(x,y)
\leftrightarrow x+yi$.  The bottom row of~\eqref{choleskygamma} is
$\gamma = \frac{b + \sqrt{D}}{2a}$, an element of the imaginary
quadratic field $K = \Q(\sqrt{D})$.  The lattice $L = L_0\,g$ has
distinguished basis $\{1, \gamma\}$.
%%% \cite[p.~149]{BoSh}.

Quadratic forms $A$, $A'$ are \emph{equivalent} if $M A M^T = A'$ for
some $M\in\GL_2(\Z)$.  They are \emph{properly equivalent} if we may
take $M \in\SL_2(\Z)$.
%% Proper equivalence is the more important notion---see
%% \cite[\S2.4]{Cox} for the history.

\subsection{Proper ideals and orders} \label{subsec:orders}

A \emph{fractional ideal} in~$K$ is a finitely generated
$\Z$-submodule of~$K$ that spans~$K$ over~$\Q$.  Fractional ideals are
lattices.
%% and lattices are fractional ideals if they are contained in some
%% quadratic subfield~$K$.
Fractional ideals $L_1, L_2$ in~$K$ are \emph{similar}, denoted $L_1
\sim L_2$, if $L_2 = \alpha L_1$ for some $\alpha \in K$, $\alpha\ne
0$.  Similar fractional ideals correspond to the same quadratic form
in~$X$.

An \emph{order}~$\O$ in a number field is a subring which contains~$1$
and which spans the number field over~$\Q$.  It is a lattice.  $K$ has
a unique \emph{maximal order} $\O_K$, its ring of integers.

The \emph{coefficient ring} of a fractional ideal~$L$ of~$K$ is $\O_L
= \{\alpha\in K \mid \alpha L \subseteq L\}$.  It is an order of~$K$.
Similar fractional ideals have the same coefficient ring.  $L$ is
similar to an ideal~$L' \subseteq \O_L$.

In quadratic fields, the orders have a particularly simple
classification \cite[\S2.7]{BoSh}. Let~$D$ be an integer
%% (positive or negative in this paragraph)
which is not a square, and let~$d$ be the squarefree part of~$D$.  Let
$K = \Q(\sqrt{D}) = \Q(\sqrt{d})$.  Let $\beta = \sqrt{d}$ if $d\equiv
2,3\bmod 4$ and $\beta = \frac{1+\sqrt{d}}{2}$ if $d\equiv 1\bmod 4$.
Then $\O_K = \Z[\beta]$.  The discriminant of $\O_K$ is called the
\emph{discriminant of}~$K$, or the \emph{fundamental discriminant}
$d_0$ for~$d$.  We have $d_0 = 4d$ if $d\equiv 2,3\bmod 4$ and $d_0 =
d$ if $d\equiv 1 \bmod 4$.  The orders of~$K$ are the rings $\O_f =
\Z[f\beta]$ for positive integers~$f$; here~$f$ is called the
\emph{conductor}.  The discriminant of $\O_f$ is $f^2 d_0$.
Conversely, to each integer~$D$, there is at most one order of
discriminant~$D$.

As in~\eqref{choleskygamma}, the primitive form $ax^2 + bxy + cy^2$
corresponds to the lattice~$L$ with basis $\{1, \gamma\}$ with $\gamma
= \frac{b + \sqrt{D}}{2a}$.  Here $\O_L$ has discriminant~$D$
\cite[p.~136]{BoSh}.

An ideal $L \subseteq \O$ may have a coefficient ring strictly larger
than~$\O$.  We say~$L$ is a \emph{proper ideal} of~$\O$ if it is an
ideal of~$\O$ and its coefficient ring equals~$\O$ \cite[p.~121]{Cox}.
A fractional ideal similar to a proper ideal is also called proper.

Let~$L$ be a proper fractional ideal of~$\O$.  Consider a matrix~$B$
whose rows are a basis of~$L$ expressed in coordinates with respect to
a basis of~$\O$.  Then $|\det B|$ does not depend on the choice of
bases and is called the \emph{norm} $N(L)$ of~$L$.  If $L\subseteq\O$,
then $N(L) = [\O : L]$.  The norm is multiplicative: $N(L_1 L_2) =
N(L_1) N(L_2)$.

The product of proper fractional ideals $L_1$, $L_2$ of~$\O$ is a
proper fractional ideal for the same~$\O$ \cite[\S2.7.4]{BoSh}.
%% $\O$ is an identity element for this product.
If $\bar L$ is the image of~$L$ under complex conjugation, then $L
\bar L = N(L)\O$, which means $\bar L$ is a multiplicative inverse
for~$L$ up to similarity.  This shows the similarity classes of proper
fractional ideals for~$\O$ form an abelian group under multiplication,
the \emph{(proper) ideal class group} of~$\O$.  This group is finite
\cite[p.~139]{BoSh}.

Complex conjugation defines an involution $\mathcal{C} \mapsto
\bar{\mathcal{C}}$ on ideal classes.  An ideal class~$\mathcal{C}$ is
\emph{ambiguous}\footnote{``Gauss, writing in Latin, used the word
\emph{anceps}, meaning either `two-headed' or `ambiguous', and with
deliberate reference to the Roman god Janus.''  \cite[p.~7]{Bu}} if
$\mathcal{C} = \bar{\mathcal{C}}$, or, equivalently, $\mathcal{C}^2 =
(1)$.

\subsection{Well-rounded ideal classes}

If the fractional ideal~$L$ is a well-rounded lattice, then any
similar fractional ideal is well rounded.  Well-roundedness depends
only on the ideal class.

\begin{proposition}
Let $\O$ be an order in an imaginary quadratic number field~$K$.  Let
$L$ be a proper fractional ideal of $\O$ that is well rounded.  Then
$L^2$ is principal.
\end{proposition}

\begin{proof}
  Let $\{v_1, v_2\}$ be a well-rounded $\Z$-basis of $L$.  Let $w =
  v_1 + v_2$.  Let $V$ be the fractional ideal $V = w^{-1} L$.
  Clearly $V^2$ is principal if and only if $L^2$ is.  By definition,
  $v_1$ and $v_2$ have the same length in~$\C$.  The generators $v_1 /
  w$ and $v_2 / w$ of $V$ have equal length, and they make equal and
  opposite angles with the real axis.  Therefore $V = \bar{V}$.  Thus
  $V^2 \sim V\bar{V} = (N(V))$ is principal.
\end{proof}

\begin{corollary} \label{finDeSiecle}
  Well-rounded classes are ambiguous.
\end{corollary}

As for the converse, ambiguous classes are not well rounded in
general.  See Section~\ref{subsec:ex1155}.

\subsection{Reduced forms}

A primitive form $ax^2 + bxy + cy^2$ is \emph{reduced}
\cite[p.~149]{BoSh} if\footnote{Other authors
%% including Lagrange, Gauss, \cite[(2.7)]{Cox}, and \cite{Bu},
choose $b \geqslant 0$ in the boundary cases $|b|=a$ or $a=c$.}
\[
|b| \leqslant a, \mathrm{\ and\ }
b \leqslant 0 \mathrm{\ if\ either\ } |b| = a \mathrm{\ or\ } a = c.
\]
The minimal non-zero value of this reduced form is~$a$.  A well-known
reduction algorithm converts any primitive form into a unique properly
equivalent reduced form of the same discriminant.  A form is reduced
if and only if its distinguished lattice basis~\eqref{choleskygamma}
is reduced.  A major result is that~\eqref{choleskygamma} gives a
bijection between the reduced forms of discriminant~$D$ and the proper
ideal class group for the order of discriminant~$D$ \cite[\S2.7]{BoSh}
\cite[\S2.A]{Cox}.

By \cite[Lemma~3.10]{Cox}, a reduced form $(a,b,c)$ is ambiguous if
and only if
\begin{equation} \label{bqfAmbig}
  b = 0, \mathrm{\ } b = -a, \mathrm{\ or\ } a = c.
\end{equation}
For reduced lattice bases $\{1,w\}$ corresponding to these three
cases, $w$ lies in Figure~\ref{fig:sl2zFundDom} as follows:
\begin{itemize}[noitemsep]
\item on the $y$-axis in the center of the fundamental domain,
\item on a vertical edge of the fundamental domain, or
\item on the curved bottom edge of the fundamental domain.
\end{itemize}

We call a primitive form $(a,b,c)$ \emph{well rounded} if the
corresponding fractional ideal is well rounded.  As we have said,
well-roundedness depends only on the class.

\begin{lemma} \label{abaIsWR}
  A reduced form $(a,b,c)$ is well rounded if and only if $a = c$.
\end{lemma}

\begin{proof}
This follows from \cite[\S2.7.7]{BoSh} (the third bullet point above).
\end{proof}

Since a homothety does not change well-roundedness, the statement of
Lemma~\ref{abaIsWR} holds even if $(a,b,c)$ is not primitive.

%% Let $L\subseteq \O_f$ be an ideal.  Say that~$L$ is \emph{prime to the
%%   conductor}~$f$ if $L + f\O = \O$.  All such~$L$
%% are proper to~$\O_f$.  Also, $L$ is prime to the conductor if and only
%% if its norm $N(L)$ is relatively prime to~$f$.  In every proper ideal
%% class for~$\O_f$, there are representative ideals prime to the
%% conductor, and the group of similarity classes of ideals prime to the
%% conductor is naturally isomorphic to the class group
%% \cite[7.C]{Cox}.  Fractional ideals relatively prime to the conductor
%% enjoy unique factorization into prime ideals \cite[Thm.~10.19]{Cohn}.

%%%%%%%%%%%%%%%%%%%%%%%%%%%%%%%%%%%%%%%%

\section{Duality}

Let $L$ be a lattice with inner product $\langle\,,\,\rangle$.  The
\emph{dual} of~$L$ is $L^* = \{w\in\R^2 \mid \langle v, w \rangle \in
\Z$ for all $v\in L\}$.  If~$L$ has $\Z$-basis $\{v_1, v_2\}$,
then~$L^*$ is a lattice with $\Z$-basis the \emph{dual basis}
$\{v_1^*, v_2^*\}$ characterized by $\langle v_i, v_j^* \rangle =
\delta_{ij}$ (Kronecker delta).  The double dual $L^{{*}{*}}$ is~$L$,
canonically.

\begin{lemma} \label{sSpBasesSameShape}
  A plane lattice and its dual have $\Z$-bases which differ from each
  other by a rotation and a homothety.
\end{lemma}

\begin{proof}
  Let $\{v_1,v_2\}$ be any basis of the lattice.  Say $\|v_1\|
  \leqslant \|v_2\|$.  Up to rotation and homothety, $v_1 = (1,0)$ and
  $v_2 = (x,y)$ with $y>0$ and $x^2 + y^2 \geqslant 1$.  The dual
  basis is $\{(1, -x/y)$, $(0, 1/y)\}$.  Calculation shows the second
  dual basis vector is shorter than the first.  To preserve
  orientation, interchange the dual vectors and multiply the longer
  one by~$-1$, obtaining the basis $w_1=(0, 1/y)$, $w_2=(-1, x/y)$
  of~$L^*$. Identify~$\R^2$ with~$\C$ as before.  The ratio $v_2 /
  v_1$ is $x+yi$, and the ratio $w_2 / w_1 = (-1 + (x/y)i)/(i/y) =
  x+yi$.  Since the ratios are equal, $\{v_1,v_2\}$ and $\{w_1,w_2\}$
  differ only by a rotation and homothety.
\end{proof}

\begin{corollary} \label{redDual}
  The dual basis of a reduced lattice basis is a reduced basis.
\end{corollary}

%% \begin{proof} 
%%   As in Section~\ref{subsec:redBases}, a reduced basis up to rotation
%%   and homothety has the form $(1,0)$, $(x,y)$ where $y>0$, $x^2 + y^2
%%   \geqslant 1$, and the length $|x|$ of the orthogonal projection of
%%   the second vector onto the first has cosine $\leqslant\frac12$.  The
%%   basis of the dual found in Lemma~\ref{sSpBasesSameShape} inherits
%%   all these properties.
%% \end{proof}
  
\begin{corollary} \label{sSpDual}
  $L$ and $L^*$ have the same number of pairs of minimal vectors.
\end{corollary}

If $L = L_0\,g$, then $L^* = L_0\,g^*$ for
$g^* = (g^{-1})^T$.
%% since the columns of $g^{-1}$ are the dual basis for the rows of $g$.

If $M \subseteq L$ are lattices, then $L^* \subseteq M^*$ and $[M^* :
  L^*] = [L : M]$.

\begin{proposition} \label{propDual}
  If $(L,M)$ is an $s$-and-$s'$ tempered perfect lattice of
  index~$\ell$ and temperament~$\tau$, then $(M^*,L^*)$ is an
  $s'$-and-$s$ tempered perfect lattice of index~$\ell$ and
  temperament $\ell/\tau$.
\end{proposition}

\begin{proof}
  By Corollary~\ref{sSpDual}, $M^*$ has~$s'$ pairs of minimal vectors
  and~$L^*$ has~$s$ pairs of minimal vectors.  A quadratic form is
  rational if and only if its dual is rational.
  Proposition~\ref{sPlusSpConverse} then implies $(M^*,L^*)$ is a
  tempered perfect form of index~$\ell$.

  We now consider the temperaments.  When $\tau=1$ (resp.,~$\ell$), the
  lattice~$L$ (resp.,~$L^*$) is hexagonal, and the result follows
  directly.  Now assume $1 < \tau < \ell$.  Let $\{v_1,v_2\}$ be a
  reduced basis of~$L$ and $\{w_1,w_2\}$ a reduced basis of~$M$.  The
  vectors $v_1,v_2,w_1,w_2$ are all distinct.  Suppose $\|v_1\|
  \leqslant \|v_2\|$ and $\|w_1\| \leqslant \|w_2\|$.  Let $v_2^\perp$
  be the orthogonal projection of $v_2$ onto the line perpendicular
  to~$v_1$.  Clearly $\|v_1\| \|v_2^\perp\| = \det L$.  Define
  $w_2^\perp$ similarly, with $\|w_1\| \|w_2^\perp\| = \det M$.  The
  proof of Lemma~\ref{sSpBasesSameShape} shows that $\|v_1^*\| =
  1/\|v_2^\perp\|$ and $\|w_1^*\| = 1/\|w_2^\perp\|$.  By definition,
  the temperament~$\tau$ for $(L,M)$ is $\tau = \|w_1\| / \|v_1\|$,
  and the temperament~$\tau^*$ for $(M^*,L^*)$ is $\tau^* = \|v_1^*\|
  / \|w_1^*\|$.  We observe
  \[
  \tau^* = \frac{\|v_1^*\|}{\|w_1^*\|} =
  \frac{\|w_2^\perp\|}{\|v_2^\perp\|} = \frac{\det M / \|w_1\|}{\det
    L / \|v_1\|} = \frac{[L : M]}{\tau}.
  \]
\end{proof}

Consider a binary quadratic form with matrix $A = g g^T$.  The
\emph{dual form} has matrix $A^* = g^* (g^*)^T$.  A calculation shows
the dual form for $ax^2 + bxy + cy^2$ is $cx^2 - bxy + ay^2$ up to
homothety.  These two forms have the same discriminant.

%%%%%%%%%%%%%%%%%%%%%%%%%%%%%%%%%%%%%%%%

\section{Doubly well-rounded forms}

A tempered perfect lattice is a pair $(L,M)$ where~$L$ and~$M$ are
fractional ideals in the same imaginary quadratic field, possibly with
different coefficient rings.  We now prove that, for doubly
well-rounded forms, $L$ and $M$ have the \emph{same} coefficient ring.
``Doubly well-rounded'' includes 3-and-3 forms as well as 2-and-2
forms.

\subsection{Coefficient ring of the lattice and sublattice}

Consider a primitive reduced form $(a,b,c)$ on~$L_0$.  Suppose this defines a
doubly well-rounded form $(L_0,M_0)$ for some sublattice~$M_0$ of
index~$\ell$.  Lemma~\ref{abaIsWR} applied to~$L_0$ shows the
form is in fact $(a,b,a)$ with $\gcd(a,b)=1$.  Let $D = b^2 - 4a^2$ be
the discriminant.  When we restrict the form to the index~$\ell$
sublattice~$M_0$, write it in coordinates $x,y$ with respect to some
$\Z$-basis of~$M_0$, and reduce it, it will have different
coefficients $(u,v,w)$, not necessarily primitive.  Since~$M_0$ is
well rounded, the form is $(u,v,u)$.  In this situation, $\tau^2 =
u/a$.

\begin{proposition} \label{oldconj6}
  Consider a doubly well-rounded reduced form, with notation as above.
  Then $\gcd(u,v) = \ell$.  Letting $u_1 = u/\ell$ and $v_1 = v/\ell$,
  the form $(u_1,v_1,u_1)$ is a primitive and reduced form on~$M_0$ of
  discriminant~$D$.  The lattice and sublattice have the same
  coefficient ring.  The temperament~$\tau$ of this doubly
  well-rounded form satisfies
  \[
  \tau^2 = \ell \frac{u_1}{a}.
  \]
\end{proposition}

\begin{proof}
  Temporarily change $(a,b,a)$ to a properly equivalent form
  $(a_1,b_1,c_1)$ with $\ell\nmid a_1$; this is possible by
  \cite[(2.25)]{Cox}.  The matrix of the form on~$M_0$ is
  $H \begin{bmatrix} a_1 & b_1/2 \\ b_1/2 & c_1 \end{bmatrix} H^T$,
  where~$H$ is one of the $\ell+1$ matrices in~\eqref{heckeMatsEll}.
  \cite[Prop.~7.3]{Bu} addresses this situation.  It shows that the
  form on~$M_0$ has discriminant either~$D$ or~$\ell^2 D$.
  Furthermore, it shows that the discriminant is~$D$ if and only if
  the form on~$M_0$ is primitive, and that in the imprimitive case the
  form is~$\ell$ times a primitive form of discriminant~$D$.  Thus to
  prove our proposition, it suffices to show that the form on~$M_0$
  does not have coefficient ring of discriminant $\ell^2 D$.

  Consider the chain of inclusions $\ell L_0 \subset M_0 \subset L_0$.
  Both inclusions have index~$\ell$, so the result from~\cite{Bu}
  applies to each inclusion.  The discriminant of the coefficient ring
  for~$M_0$ is equal to either~1 or~$\ell^2$ times the discriminant of
  the coefficient ring for~$L_0$, and the discriminant of the
  coefficient ring for~$\ell L_0$ is equal to either~1 or~$\ell^2$
  times the discriminant of the coefficient ring for~$M_0$.  But the
  coefficient rings for $L_0$ and $\ell L_0$ are equal, since the
  lattices are homothetic.  Thus all three coefficient rings are the
  same.
\end{proof}

\subsection{Characterization of doubly well-rounded forms by ideals}

Since well-roundedness depends only on the ideal class, we would like
to use ideal classes to compare the two lattices in a doubly
well-rounded form.

\begin{proposition} \label{frakpellbetween}
Given a doubly well-rounded form, if~$\O$ is the common coefficient ring of
its lattices $(L,M)$, then~$\O$ contains a prime proper ideal
$\p_\ell$ of norm~$\ell$, and $M = \p_\ell L$.
\end{proposition}

\begin{proof}
  The multiplicative inverse of~$L$ is proper and equals
  $\frac{1}{N(L)} \bar L$.  Thus $\p_\ell = \frac{1}{N(L)} \bar{L} M$
  is a proper fractional ideal of~$\O$ satisfying $\p_\ell L = M$.
  Since $M\subset L$, we have $\p_\ell \subseteq\O$, and we may drop
  the word ``fractional'' for~$\p_\ell$.  Since $[L:M] = \ell$, we
  have $N(\p_\ell) = \ell$.  The ideal~$\p_\ell$ is prime, indeed
  maximal, because $|\O/\p_\ell|$ is prime.
\end{proof}

Combining Propositions~\ref{oldconj6} and~\ref{frakpellbetween}, we
have characterized doubly well-rounded forms in terms of the class
group, as follows.

\begin{corollary} \label{oldconj5}
  Let $\mathfrak{a}$ be an ideal of $\O$ in a well-rounded class
  $\mathcal{C}$. Suppose~$\p_\ell$ is a prime proper ideal of $\O$ of
  norm $\ell$, and suppose $\p_\ell\mathfrak{a}$ is in a well-rounded
  class $\mathcal{C}'$. Then $(L,M) = (\mathfrak{a},
  \p_\ell\mathfrak{a})$ is doubly well rounded.  Conversely, every
  doubly well-rounded form occurs in this way.
\end{corollary}

%%%%%%%%%%%%%%%%%%%%%%%%%%%%%%%%%%%%%%%%

\section{Tempered perfect forms from the Eisenstein integers}

Because of duality, it suffices to consider $s$-and-$s'$ forms where
$s \geqslant s'$. By Corollary~\ref{sPlusSp}, $s + s' \geqslant 4$. In
this section we characterize the temperaments~$\tau$ of 3-and-$s'$
forms.  The techniques differ from those for 2-and-2 forms, since for
3-and-$s'$ forms we work exclusively in the hexagonal lattice, the
Eisenstein integers $\Z[\omega]$.  We will return to 2-and-2 forms in
Section~\ref{sec:22asFuncEll}.

We always identify vectors $v\in\R^2$ with the corresponding $z\in\C$.
The \emph{norm} $N(z)$ of $z\in\C$ is its squared length $\|v\|^2 =
\langle v,v \rangle = |z|^2$.  When $z\in\Z[\omega]$ (or, in the later
sections, in other orders~$\O$), we also identify~$z$ with its
principal ideal~$(z)$.  The norm $N(z)$ equals the ideal norm
of~$(z)$.

For 2-and-2 forms, we saw in Proposition~\ref{oldconj6} that~$L$
and~$M$ are ideals for a common order~$\O$.  For 3-and-1 forms,
however, any sublattice~$M$ of $\Z[\omega]$, not necessarily an ideal,
with minimal norm~$N$ in an appropriate range will yield a tempered
perfect form with $\tau^2 = N$.  Classifying 3-and-1 forms is thus a
matter of characterizing the minimal norms that arise for sublattices
$M \subset \Z[\omega]$ of prime index $\ell$.

\subsection{The Eisenstein integers} \label{subsec:defEis}

In $\Z[\omega]$, the three pairs of minimal vectors are $\pm 1, \pm
\omega, \pm \omega^2$.  These six are also the units of $\Z[\omega]$.
%% The fundamental parallelogram of this lattice has area $\sqrt{3}/2$.
The norm $N(a+b\omega) = a^2 - ab + b^2$.

Let~$\p_p$ be a prime ideal of $\Z[\omega]$ over a rational prime~$p$.
Write $\p_p = (\pi_p)$ with $\pi_p \in \Z[\omega]$.  One checks that
$\p_3 = (1-\omega)$ is ramified, meaning $(3) = \p_3^2$.  The primes
$p \equiv 1 \bmod 6$ split: $(p) = \p_p\overline{\p_p} =
(\pi_p)(\overline{\pi_p})$, where $\p_p \ne \overline{\p_p}$ and
$N(\pi_p) = p$.  The primes $p \equiv 5 \bmod 6$ remain inert, meaning
$\p_p = (p)$.

We call $x = a + b\omega \in \Z[\omega]$ a \emph{primitive
representation of}~$n$ if $gcd(a,b) = 1$ and $N(x) = n$.

\begin{proposition} \label{numRepsEis}
Suppose $n$ is primitively represented by $\Z[\omega]$. Let $S$ be the
set of rational primes $p \equiv 1 \bmod 6$ such that $p|n$. Then the
prime factorization of $n$ has the form $3^r \prod_{p \in S} p^{r_p}$
with $r=0$ or~$1$, and there are $6 \cdot 2^{|S|}$ primitive
representations of~$n$.
\end{proposition}

\begin{proof} \cite[p.~176]{NZM}. \end{proof}

\subsection{Well-rounded sublattices of the Eisenstein integers}

\begin{proposition} \label{noNonPrinceWR}
All well-rounded sublattices~$M$ of $\Z[\omega]$ with prime
index~$\ell$ are ideals of $\Z[\omega]$.
\end{proposition}

\begin{proof}
  This is the $L=\Z[\omega]$ case of Proposition~\ref{oldconj6}.
\end{proof}

\begin{corollary}
There are no $3$-and-$2$ forms and no $2$-and-$3$ forms.
\end{corollary}

\begin{proof}
  If $(L,M)$ is a 3-and-$s'$ form for $s' \geqslant 2$, then, by
  Proposition~\ref{noNonPrinceWR}, $M$ is an ideal of $\Z[\omega]$.
  Ideals are closed under multiplication by the six units, so $s'=3$.
  For 2-and-3 forms, apply Proposition~\ref{propDual}.
\end{proof}

\begin{corollary} \label{all33}
  All $3$-and-$3$ forms are equivalent to $(L,M) = (\Z[\omega],
  \p_\ell)$ where $\ell=3$ or $\ell \equiv 1 \bmod 6$.
\end{corollary}

\begin{proof}
  Since $\Z[\omega]$ is a principal ideal domain, it has exactly one
  ideal class.  This class is well rounded.  The rest follows from
  Corollary~\ref{oldconj5}.
\end{proof}

\subsubsection{Example of a 3-and-3 form}
\begin{figure}
  \begin{center}
    \includegraphics[scale=0.2]{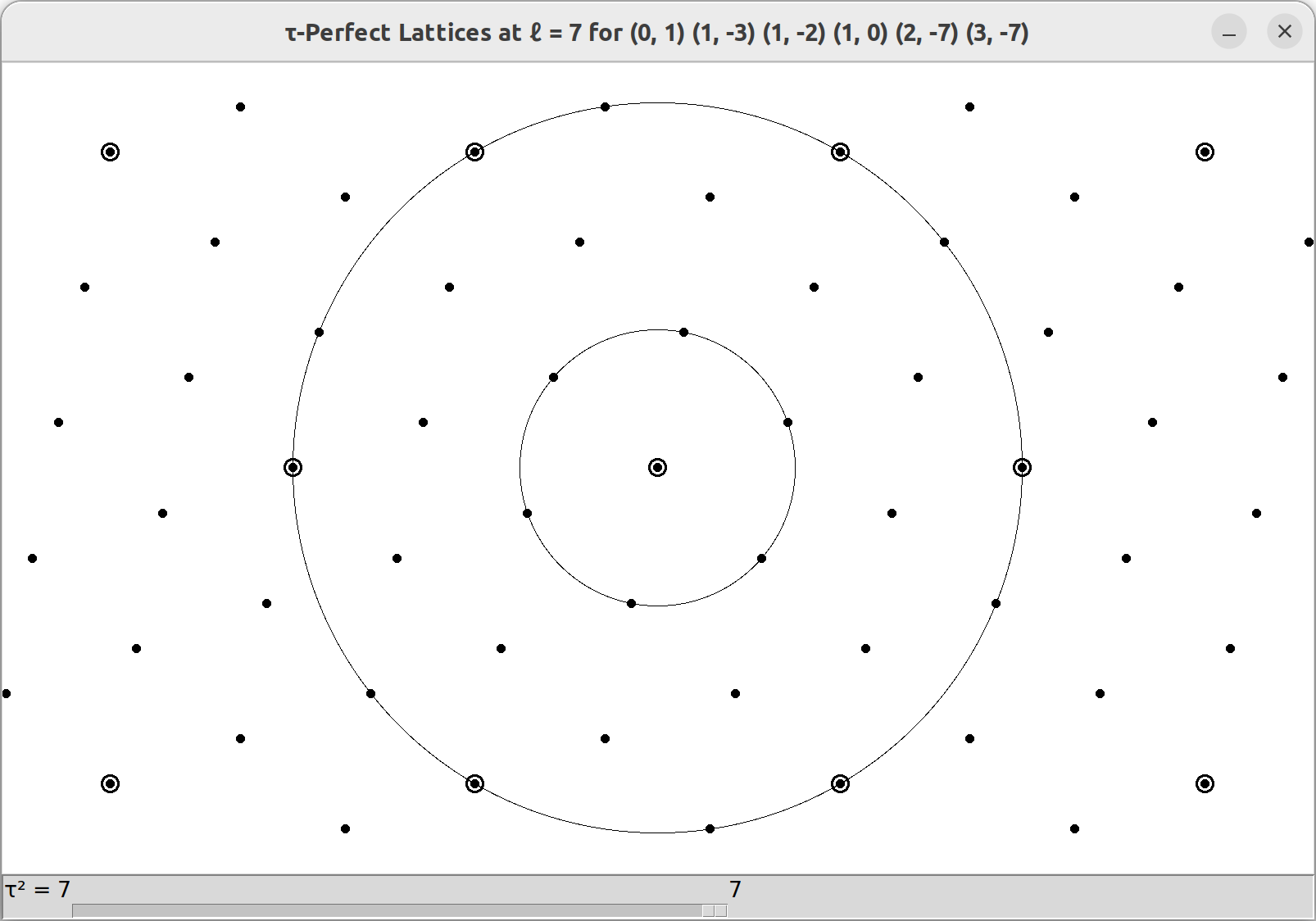}
  \end{center}
  \caption{The 3-and-3 form from the prime ideal over~7 in
    $\Z[\omega]$.}
  \label{fig:exEis7}
\end{figure}

In Figure~\ref{fig:exEis7}, the dots are $L = \Z[\omega]$.  The
sublattice~$M$ of~$\odot$ dots is the prime ideal $\p_7 = (2-\omega)$.
Here $\ell = 7$ and $\tau = \sqrt7$.

\subsection{General sublattices of the Eisenstein integers}

\begin{lemma} \label{MinkowskiBasisMaker}
  Suppose an index~$\ell$ sublattice $M$ of $\Z[\omega]$ has linearly
  independent vectors $v_1$, $v_2$ representing positive integers
  $n_1$, $n_2$ both less than or equal to~$\ell$. Then $\{v_1, v_2\}$
  is a reduced basis of~$M$.
\end{lemma}

\begin{proof}
Let $\arg(v,w)$ denote the angle between $v$ and $w$, measured
counterclockwise from $v$.  The area of the fundamental parallelogram
for $\Z[\omega]$ is $\sqrt{3}/2$.  The pair $\{v_1, v_2\}$ must
generate a sublattice~$M'$ of index at least $\ell$, hence of area at
least $(\sqrt{3}/2)\ell$.  Thus $\|v_1\| \|v_2\| |\sin\arg(v_1,v_2)|
\geqslant (\sqrt{3}/2)\ell$ with $\|v_i\| \leqslant \sqrt{\ell}$,
whence $\arg(\pm v_1, \pm v_2) \geqslant \pi/3$. We may negate $v_1$
or $v_2$ if necessary to suppose $\pi/3 \leqslant \arg(v_1, v_2)
\leqslant \pi/2$. Suppose $\{v_1, v_2\}$ as above do not form a
reduced basis, so that there exists $v_3$ with $\|v_3\| < \max(n_1,
n_2)$. Then if we consider the order of $v_1, v_2, v_3$ sweeping
counterclockwise from $v_1$, we have two cases:
\begin{gather*}
  \arg(v_1, v_2) = \arg(v_1, v_3) + \arg(v_3, v_2) \leqslant \pi/2 \\
  \arg(v_2, -v_1) = \arg(v_2, v_3) + \arg(v_3, -v_1) \leqslant 2\pi/3
\end{gather*}
In the first case either $\arg(v_1,v_3) \leqslant \pi/4$ or
$\arg(v_3,v_2) \leqslant \pi/4$, a contradiction, since this would
imply they generate a sublattice of index smaller than $\ell$. In the
second case, either $\arg(v_2, v_3) \leqslant \pi/3$ or $\arg(v_3,
-v_1) \leqslant \pi/3$, which is only possible if the sublattice is
the principal ideal $(v_1)$; in this case $\{v_1, v_2\}$ would still
be a reduced basis, a contradiction.
\end{proof}

\begin{corollary} \label{lemma2}
Let $v$ be a primitive vector in $\Z[\omega]$, and suppose $\|v\|^2 <
\ell$. Then $v$ is contained in exactly one index~$\ell$ sublattice
$M$, and $v$ is one of the vectors in a reduced basis of~$M$.
\end{corollary}

\begin{proof}
Lemma~\ref{uniqueSublatt} gives the existence and uniqueness of~$M$.
Since~$v$ is primitive, we can complete it to some lattice basis
$\{v,w\}$ of~$M$, where~$w$ is chosen to be as short as possible.  If
$\|w\| \geqslant\ell$, then $\|v\| < \|w\|$ and the basis is reduced.
If $\|w\| < \ell$, the result follows from
Lemma~\ref{MinkowskiBasisMaker}.
\end{proof}

%% As always, the reduced basis is unique up to sign changes and
%% interchanging the vectors.  The only exception is the hexagonal
%% lattice, with its extra six-fold symmetry.
 
\begin{definition}
$E_\ell(n)$ is the set of all index~$\ell$ sublattices of $\Z[\omega]$
  that contain a primitive representation of~$n$.
\end{definition}

\begin{corollary}
Let $n < \ell$ be primitively representable. Let $S$ be the set of
rational primes $p \equiv 1 \bmod 6$ such that $p|n$. Then
$|E_\ell(n)| = 3 \cdot 2^{|S|}$.
\end{corollary}

\begin{proof}
Let $x_1$ and $x_2$ be primitive representations of $n$ that are
distinct up to units.  By Lemma~\ref{MinkowskiBasisMaker}, they form a
reduced basis of the sublattice they generate.  This basis is well
rounded.  By Proposition \ref{noNonPrinceWR}, however, being distinct
up to units, they cannot lie in a common index~$\ell$ sublattice.  We
conclude there is an injective map from primitive representations~$x$
of~$n$ (up to units) to sublattices~$M$ of index~$\ell$.  The number
of primitive representations~$x$ comes from
Proposition~\ref{numRepsEis}.  Since $n$ is either the first or second
shortest norm in a given~$M$, we see that for any choice of
representation, rotation by roots of unity (except for $\pm 1$)
generates new sublattices. Thus in total there are $3 \cdot 2^{|S|}$
distinct sublattices representing $n$.
\end{proof}

The next results will help determine which $n$ arise as minimal norms
represented by index $\ell$ sublattices of $\Z[\omega]$.

\begin{lemma} \label{noPrimeSharing}
Suppose $\{v,w\}$ is a reduced basis of an index~$\ell$ sublattice~$M$
of $\Z[\omega]$ and at least one of the principal ideals $(v)$, $(w)$
in $\Z[\omega]$ is not divisible by a prime $\p_\ell$
over~$\ell$. Then~$(v)$ and~$(w)$ share no prime ideals in their
factorization.
\end{lemma}

\begin{proof}
  Suppose $(v) = \p(v')$ and $(w) = \p(w')$ where~$\p = (u)$ is a
  prime ideal of $\Z[\omega]$ over~$p$ and $v', w' \in\Z[\omega]$.
  Since~$v$ and~$w$ are primitive, $\p$ is not inert.  Then $\{\bar{u}
  v, \bar{u} w\}$ generates a sublattice of $\Z[\omega]$ of index
  $p\ell$.  Hence $\{v',w'\}$ generates a sublattice of $\Z[\omega]$
  of index $\ell/p$.  This implies $p|\ell$, contradicting the
  hypothesis.
\end{proof}

\begin{proposition} \label{prop:e06}
Suppose $n_1$, $n_2$ are distinct and both $n_1, n_2 < \ell$. Then
$|E_\ell(n_1) \cap E_\ell(n_2)| = 0$ or $6$.
\end{proposition}

\begin{proof}
Suppose~$M_1$ and~$M_2$ are both in $E_\ell(n_1) \cap E_\ell(n_2)$.
By Lemma~\ref{MinkowskiBasisMaker}, they have bases $\{v_1,w_1\}$ and
$\{v_2,w_2\}$ where $\|v_1\| = \|v_2\| = n_1$ and $\|w_1\| = \|w_2\| =
n_2$.  The principal ideals $(v_1)$ and $(v_2)$ differ in their prime
ideal factorization only by taking conjugates of certain ideals. If we
define $S_v$ as the set of prime ideals at which $(v_1), (v_2)$
differ, and $S_w$ analogously for $(w_1), (w_2)$, we may write
(relabeling some $\p_i$ as $\overline{\p_i}$ and vice versa if
necessary)
\[
(n_1 n_2^{-1}) =
(v_1)(v_2)^{-1} =
\prod_{s \in S_v} (\overline{\p_s} \p_s^{-1})^{r_s} =
\prod_{s \in S_w} (\overline{\p_s} \p_s^{-1})^{r_s} =
(w_1)(w_2)^{-1}.
\]
By Lemma~\ref{noPrimeSharing}, both $\{\p_s \mid s \in S_v\}$ and
$\{\overline{\p_s} \mid s \in S_v\}$ must be disjoint from $\{\p_s
\mid s \in S_w\} \cup \{\overline{\p_s} \mid s \in S_w\}$. Thus if the
two products above are equal, it must be the case that $S_v$ and $S_w$
are empty.  Thus in fact $M_2$ differs from $M_1$ only by a unit.

Starting with any sublattice, three sublattices are obtained via
multiplication by units, and three more obtained via complex
conjugation followed by multiplication by units. Thus $|E_\ell(n_1)
\cap E_\ell(n_2)| = 6$ when it is not~$0$.
\end{proof}

\subsubsection{Example of a 3-and-1 form}
\begin{figure}
  \begin{center}
    \includegraphics[scale=0.2]{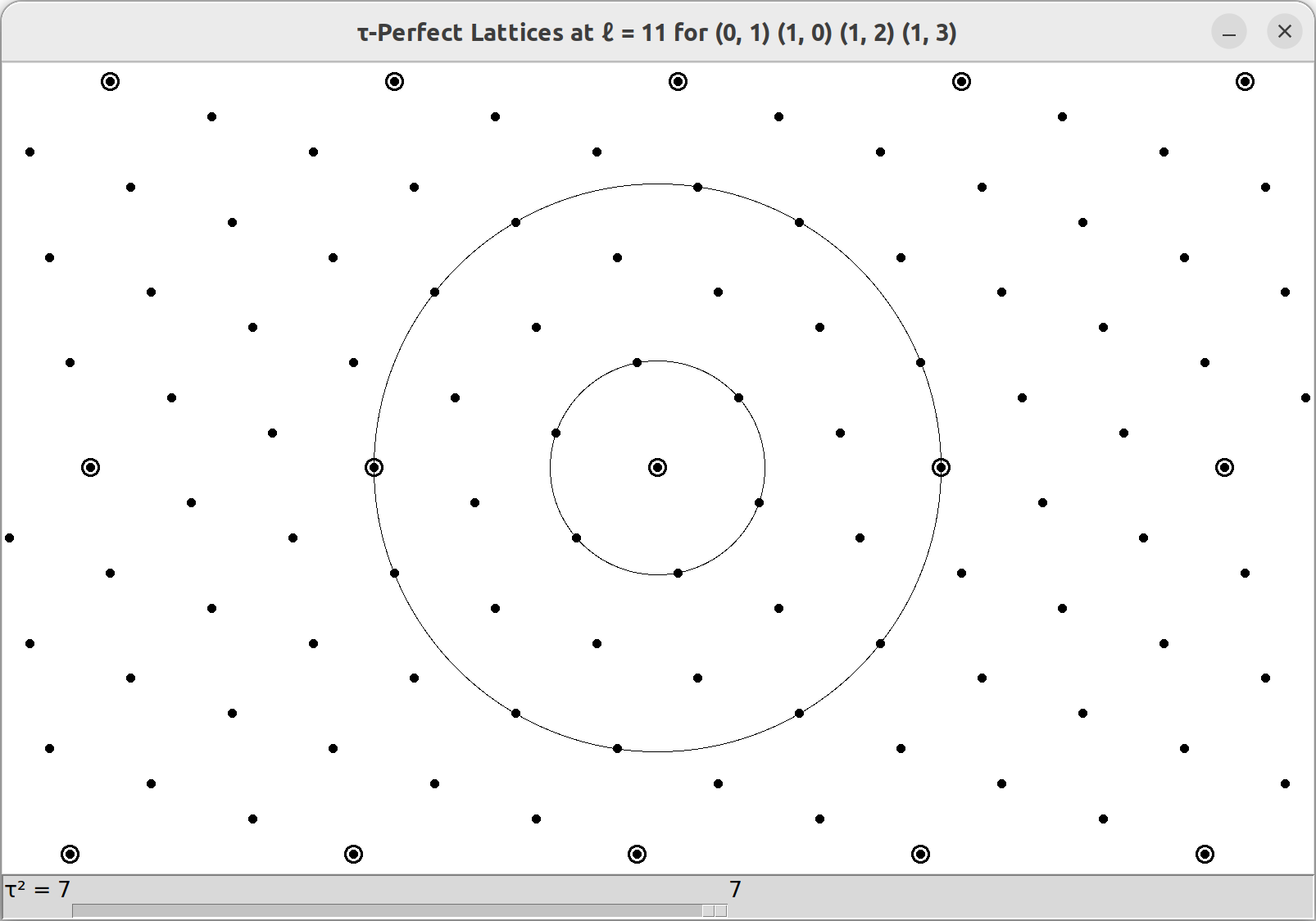}
  \end{center}
  \caption{A 3-and-1 form.  The~$\odot$ sublattice is not an ideal of
    $\Z[\omega]$.}
  \label{fig:exEis11}
\end{figure}

In Figure~\ref{fig:exEis11}, the dots are $L = \Z[\omega]$.  The
sublattice~$M$ of~$\odot$ dots has one pair of minimal vectors
$\pm(2-\omega)$, and $\tau = \sqrt7$, just as in
Figure~\ref{fig:exEis7}.  However, all other vectors in~$M$ are
longer; no other pair is on the outer circle.  $M$ is not an ideal of
$\Z[\omega]$.  The index $\ell = 11$.  Since $7<11$,
Corollary~\ref{lemma2} guarantees that $v = 2-\omega$ is contained in
exactly one index~11 sublattice and is the shortest vector in a
reduced basis $\{v,w\}$ of that sublattice.
Proposition~\ref{prop:e06} says that if we multiply~$M$ by powers
of~$\omega$, and take complex conjugates, then we get six distinct
lattices in $\Z[\omega]$ with the same properties.

\subsection{Algorithm for 3-and-1 temperaments}

Fix a prime $\ell$, and suppose that one would like to compute all
temperaments~$\tau$ that correspond to $3$-and-$1$ forms.  $\tau^2 \in
[1, \ell^2]$ by Section~\ref{bdryTau}.  However, we only need to
consider half that range:

\begin{proposition} \label{noTemps}
 All $3$-and-$1$ forms have $\tau^2 < \ell$.
\end{proposition}

\begin{proof}
A 3-and-1 form with temperament $\tau$ can only exist when $L =
\Z[\omega]$ has an index~$\ell$ sublattice $M$ with reduced basis
$\{v,w\}$ with $\|v\| < \|w\|$ and $\|v\|^2 = \tau^2$.  Since
we may assume the angle~$\theta$ between $v$ and $w$ satisfies $\pi/3
\leqslant \theta \leqslant 2\pi/3$, and because $\|v\| <
\|w\|$, the area of the fundamental parallelogram for~$M$ is strictly greater
than the area of the fundamental parallelogram for the
principal ideal~$(v)$. Therefore $\ell = [\Z[\omega] : M] >
[\Z[\omega] : (v)] = \|v\|^2$.
\end{proof}

\begin{corollary}
  All $1$-and-$3$ forms have $\tau^2 > \ell$.
\end{corollary}

\begin{proof}
  Apply Proposition~\ref{propDual}.
\end{proof}

The case of $1 \leqslant \tau^2 < (\sqrt{3}/2)\ell$ also has a simple
answer, but we require a supporting lemma.

\begin{lemma} \label{lemma1}
  Let $M$ be an index~$\ell$ sublattice of the Eisenstein lattice.
  Let~$v$ be a vector of minimal length within~$M$, and let~$w$ be a
  vector of minimal length in $L - \mathrm{span}(v)$.

  (a) If $\|v\|^2 \leqslant (\sqrt{3}/2)\ell$, then $\|w\|^2 \geqslant
  (\sqrt{3}/2)\ell$.

  (b) If $\|v\|^2 \leqslant (3/4)\ell$, then $\|w\|^2 \geqslant \ell$.
\end{lemma}

\begin{proof}
We derived the formula for the area of the parallelogram spanned by
$v$ and $w$ during the proof of Lemma~\ref{MinkowskiBasisMaker}:
$\|v\| \|w\| |\sin\theta| = \ell \frac{\sqrt{3}}{2}$.  Suppose
$\|v\|^2 = \mu\ell$.  Then $\|w\|^2 = \frac{3\ell^2}{4\|v\|^2
  \sin^2\theta} \geqslant \frac{3\ell^2}{4 (\mu\ell)(1)} =
\frac{(3/4)\ell}{\mu}$.  Both parts follow.
\end{proof}

\begin{theorem} \label{thm:eis}
Let $\ell$ be a prime.  For each~$\beta\in\Z$ with $\beta <
(\sqrt{3}/2)\ell$ which is primitively represented by an element of
$\Z[\omega]$, there is a 3-and-1 form with $\tau^2 = \beta$.
\end{theorem}

\begin{proof}
  Take $v\in\Z[\omega]$ with $N(v) = \beta$.  By
  Lemma~\ref{uniqueSublatt}, there is a unique index~$\ell$
  sublattice~$M$ containing~$v$.  By Lemma~\ref{lemma1}~(a), $v$ is
  (strictly) a minimal vector of~$M$.  By
  Proposition~\ref{sPlusSpConverse}, this gives a 3-and-1 form.
\end{proof}

We have now characterized the temperaments coming from 3-and-1 forms
except in the range $\tau^2 \in [(\sqrt{3}/2)\ell, \ell]$.
Algorithm~\ref{alg:eis} below will give us all the~$\tau^2$ in this
range.

Before discussing the algorithm, we need a supporting proposition.
Let $M$ be an index $\ell$ sublattice with minimal vector $v = a +
b\omega$.  Take~$c,d$ such that $ad - bc = 1$. Then $\begin{bmatrix} a
  & b \\ \ell c & \ell d \end{bmatrix}$ has determinant~$\ell$, which
implies by Lemma~\ref{uniqueSublatt} that the vectors $\{a + b\omega,
\ell(c + d\omega)\}$ generate $M$. Let $w = c + d\omega$. Without loss
of generality, take the angle~$\theta$ between~$v$ and~$w$ to be less
than~$\pi$.  Let $z \in \C$ be perpendicular to~$v$ and of length
$\|z\| = \frac{\sqrt{3}\ell}{2\|v\|}$.  Thus $z$ is the height of the
parallelogram spanned by $v$ and $\ell w$.  Standard arguments with
lattices give the following.

\begin{proposition}
\label{computeNextSmallest}
With $M, v, w, z$ as above, the second shortest vector of $M$ is
$\hat{w} = \ell w + tv$, where $t$ is the nearest integer to $t_0 =
\frac{h - w\ell}{v} \in \R$.
\end{proposition}

As we have said, we are now studying the range $\tau^2 \in
[(\sqrt{3}/2)\ell, \ell]$.  However, Algorithm~\ref{alg:eis} works
more naturally on the superset $[(3/4)\ell, \ell]$.  We need the
larger range because the inductive mechanism in the algorithm must
begin at $(3/4)\ell$, due to the limitations expressed in
Lemma~\ref{lemma1}.

\begin{algorithm} \label{alg:eis}
Given a prime~$\ell$, this algorithm determines all
temperaments~$\tau$ in the range $\tau^2 \in [(3/4)\ell, \ell]$ that
come from 3-and-1 forms. Let $\{\beta_\alpha\}_{\alpha = 1}^\infty$ be
the ordered list of integers which are primitively represented as
norms of Eisenstein integers, starting with $\beta_1 = 1$, $\beta_2 =
3$, $\beta_3 = 7$, \dots. Let $n$ be the last index for which
$\beta_{n} < (3/4)\ell$ and~$N$ the last index for which $\beta_N <
\ell$.
\begin{enumerate}[noitemsep]
\item $[$Create Table$]$ Initialize a list~$\mathcal{L} =
  \{\beta_{n+1}, \beta_{n+2}, ..., \beta_{N}\}$. In a table format,
  head the columns with the norms in $\mathcal{L}$, and, under each
  norm, write down all primitive vectors representing that norm, where
  vectors on the list are considered equivalent if they are related by
  a unit. In effect, each element of the table is a
  ``$\mathbb{Z}/6\mathbb{Z}$ orbit" of the action of the
  units. Initialize~$T$, the set of~$\tau^2$ which will be output, to
  the empty set.
\item $[$Compute$]$ If~$\mathcal{L}$ is empty, output~$T$ and
  terminate the algorithm.  Otherwise, choose~$v$ under the smallest
  norm~$\beta_{\alpha_0}$ available in~$\mathcal{L}$.  Using
  Proposition~\ref{computeNextSmallest}, compute~$w$ and its norm
  $\beta_{\alpha_1}$.
\item $[$Record$]$ Add $\beta_{\alpha_0}$ to~$T$.
\item $[$Cross Out$]$ Cross out the $\Z/6\Z$ orbits corresponding to
  $v$ and $w$ under the norms $\beta_{\alpha_0}$ and
  $\beta_{\alpha_1}$, respectively. If $\beta_{\alpha_1}$ is out of
  range (not on the list), only cross out the orbit corresponding to
  $v$. If, after crossing out the orbits, a column has been completely
  exhausted, delete its corresponding norm from $\mathcal{L}$. Go to
  step~2.
\end{enumerate}
\end{algorithm} 

\begin{proof}
The main hurdle is to check that every vector $v$ chosen in step~2 is
in fact the minimal vector of the unique index $\ell$ sublattice that
contains $v$. Once this is verified, the fact that $\beta_{\alpha_0}$ is a
temperament follows immediately from the definition. To show that the
algorithm successfully identifies \emph{all} of the temperaments in
the range, we also need to show that all sublattices with minimal
vectors in this range have been identified by the algorithm.

Let $v$ be the first vector chosen by step~2 of the algorithm, and let
$M$ be the unique index $\ell$ sublattice that contains $v$. Due to
Corollary~\ref{lemma2}, $v$ is contained in the reduced basis of
$M$. Let $\hat{w}$ be the other vector in this basis. Suppose that
$\|\hat{w}\| < \|v\|$. Since $v$ has norm $\beta_{n + 1}$, the norm of
$\hat{w}$ must be less than $(3/4)\ell$, which would imply by
Lemma~\ref{lemma1}~(b) that $v$ has norm greater than $\ell$. Now
suppose inductively that after $n$ cycles of steps $1$ through $4$ of
the algorithm, all vectors $v$ chosen in step~2 have been minimal
vectors. Then we wish to prove that in the $(n+1)$st cycle, step~2
produces another minimal vector. Suppose that this is not true. Let
the true minimal vector be denoted by $\hat{v}$. Then $\hat{v}$
satisfies $\|\hat{v}\| < \|v\|$ by Proposition~\ref{noNonPrinceWR}. It
also satisfies $\|\hat{v}\|^2 > (3/4)\ell$: otherwise $v$ would have
norm greater than $\ell$ by Lemma~\ref{lemma1}~(b). These bounds imply
that $\hat{v}$ was previously dealt with in the algorithm, which means
that the orbit corresponding to $v$ should have already been crossed
off.

All temperaments are reached by this process by virtue of the
correspondence between minimal vectors of index~$\ell$ sublattices and
$3$-and-$1$ forms: all minimal vectors of sublattices with norm in the
range $[(3/4)\ell,\ell]$ have been identified because of
Corollary~\ref{lemma2}.
\end{proof}

%%%%%%%%%%%%%%%%%%%%%%%%%%%%%%%%%%%%%%%%

\section{Two-and-two Forms} \label{sec:22asFuncEll}

\subsection{The Gaussian integers}

Before discussing 2-and-2 forms in general, we describe the most
common ones, those coming from the Gaussian integers $\Z[i]$.  Since
$L=\Z[i]$ is the square lattice, it is well rounded, with $s=2$.  The
ring is a principal ideal domain, so it has only one ideal class, the
well-rounded class.  By Corollary~\ref{oldconj5}, all 2-and-2 forms
for $\Z[i]$ have $M$ equal to a prime ideal~$\p_\ell$, where~$\ell$
splits or is ramified in $\Z[i]$, and where $\tau = \sqrt{\ell}$.  The
only ramified prime is $\p_2 = (1+i)$.  An odd~$\ell$ splits if and
only if $\ell\equiv 1 \bmod 4$.  Since~$M$ is an ideal, it is closed
under multiplication by~$i$, making it a square lattice in its own
right.

\begin{remark}
The only well-rounded \emph{principal} ideals come from the Eisenstein
and Gaussian integers.  All other 2-and-2 forms will come from
orders~$\O$ of class number greater than one.
\end{remark} 

\subsubsection{Example of a 2-and-2 form}

\begin{figure}
  \begin{center}
    \includegraphics[scale=0.2]{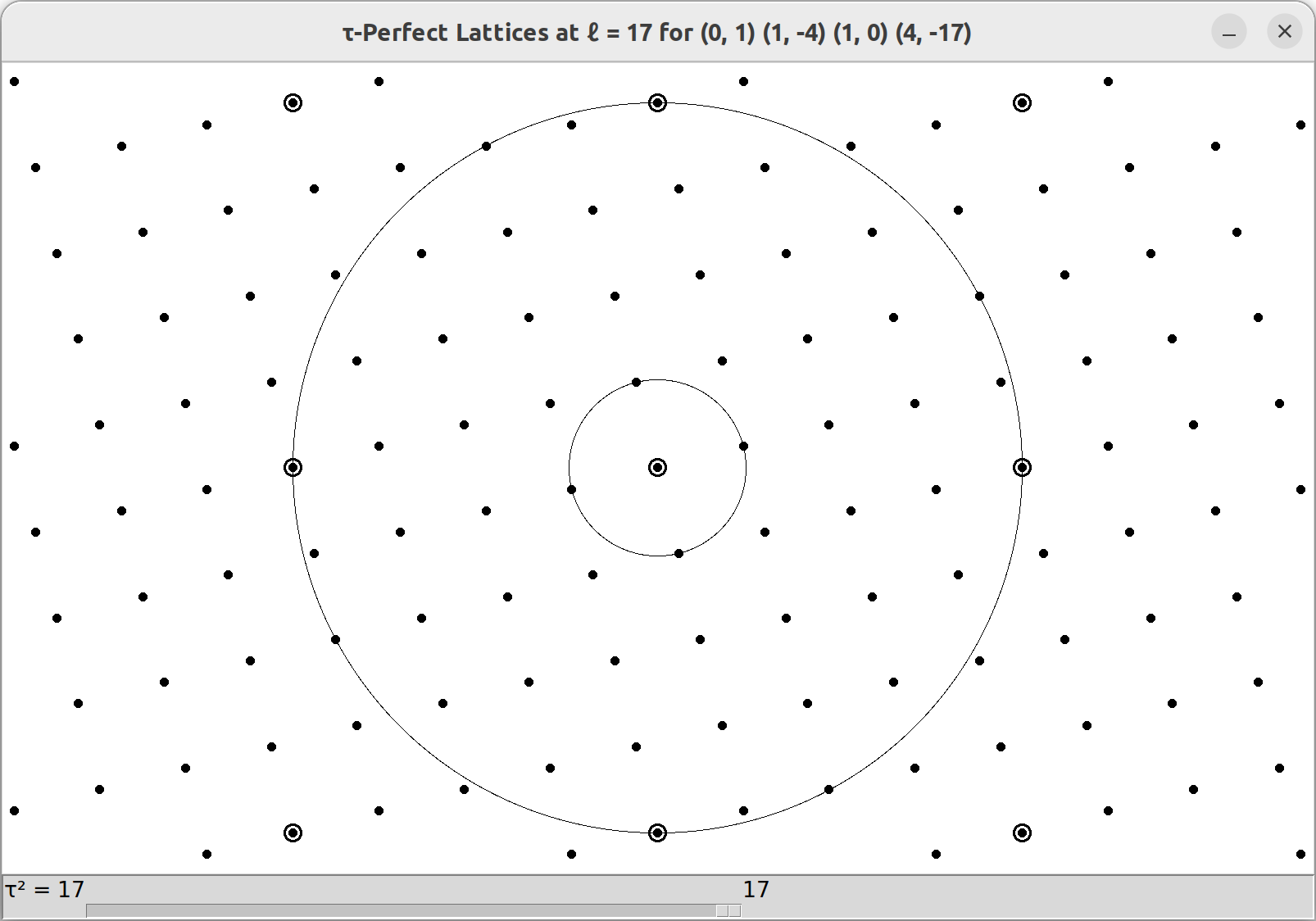}
  \end{center}
  \caption{A 2-and-2 form from $\Z[i]$: a square lattice in a square
    lattice.}
  \label{fig:exGauss17}
\end{figure}

In Figure~\ref{fig:exGauss17}, the dots are $L = \Z[i]$.  The
sublattice~$M$ of~$\odot$ dots is the principal prime ideal~$\p_{17} =
(1+4i)$.  The rational prime~17 splits as $\p_{17}\overline{\p_{17}} =
(1+4i)(1-4i)$.  We have $\ell = 17$ and $\tau = \sqrt{17}$.

\subsection{Finding all 2-and-2 forms}

%% For small~$\ell$, forms coming from the Eisenstein and Gaussian
%% integers make up all of the tempered perfect forms.  For~$\ell$ in
%% general, however, 2-and-2 forms may come from any number of imaginary
%% quadratic number fields.

Two main questions frame our work on 2-and-2 forms.  First,
given~$\ell$, we wish to show there will be only finitely many 2-and-2
forms up to $\Gamma_0(\ell)$-equivalence.  The compactness arguments
of~\cite{MM20} already show there are only finitely many, but we
prefer a number-theoretic characterization.
%% We know from Proposition~\ref{oldconj6} that each 2-and-2 form is a nested pair of proper (fractional) ideals in an order of discriminant~$D$.
Corollary~\ref{all33} answered the corresponding question for 3-and-3
forms.  Theorem~\ref{thm:eis} and Algorithm~\ref{alg:eis} answered it
for 3-and-1 forms.

The second main question is, given an order~$\O$ of discriminant~$D$,
for which~$\ell$ will a 2-and-2 form arise from~$\O$?  This will be
answered in Theorem~\ref{thm:infManyEll}.

We now answer the first main question for 2-and-2 forms.

\begin{theorem} \label{thm:finManyDForEll}
  Fix a prime~$\ell$.  Let~$D$ be the discriminant of an order~$\O$ in
  an imaginary quadratic number field such that a $2$-and-$2$ form
  $(L,M)$ of index~$\ell$ arises from a pair of ideal classes in~$\O$.
  Then
  \[
  |D| \leqslant 4\ell^2.
  \]
\end{theorem}

\begin{proof}
  Let $d < 0$ be a squarefree integer.  Let $K = \Q(\sqrt{d})$ with
  fundamental discriminant $d_0$.  Assume we have a 2-and-2 form whose
  coefficient ring is an order $\O = \O_f \subseteq \O_K$.  By
  Proposition~\ref{frakpellbetween}, some proper ideal $\p_\ell$
  of~$\O$ has norm~$\ell$.

By Corollary~\ref{finDeSiecle}, a class of modules in $\O_f$ that is
well rounded is ambiguous, that is, has order~1 or~2.
%% By the work in \cite[\S2.7]{BoSh} with the familiar fundamental
%% domain
%% It is easy to see from the geometry that the class of $\O_f$ itself is
%% never well-rounded, except if $f=1$ and $d\in\{-3,-1\}$ (the
%% Eisenstein or Gaussian integers).  Apart from these exceptions, a
%% well-rounded class thus has order exactly~2 in the class group.
By Corollary~\ref{oldconj5} applied to our 2-and-2 form, $\p_\ell$
carries a well-rounded class to another.  Thus $\p_\ell$ itself has
order~1 or~2 in the class group.  That is, either $\p_\ell$ or
$\p_\ell^2$ is principal.  Take $\alpha\in\O$ so that
$(\alpha)=\p_\ell$ or $(\alpha)=\p_\ell^2$ in the respective cases.
The norm $N(\alpha)$ is~$\ell$ or $\ell^2$, respectively.

The argument breaks further into cases.  First, suppose $N(\alpha) =
\ell$ and $d_0 \equiv 2,3 \bmod 4$.  We may write $\alpha = x +
yf\sqrt{d_0}$ for $x,y\in\Z$.  Taking norms, $\ell = x^2 + f^2 |d_0|
y^2$.  We cannot have $y=0$ because~$\ell$ is not a square.  Thus $|y|
\geqslant 1$, and we have $f^2 |d_0| = |D| \leqslant \ell$.

Next, suppose $N(\alpha) = \ell$ and $d \equiv 1 \bmod 4$.  We may
write $\alpha = x + yf\beta$ where $\beta = \frac{1+\sqrt{d_0}}{2}$.
Taking norms, $\ell = (x+yf/2)^2 + y^2f^2|d_0|/4$.  Again, $y=0$ is
impossible, and from $|y|\geqslant 1$ we deduce $|D| \leqslant 4\ell$.

Next, suppose $N(\alpha) = \ell^2$.  If $y\ne 0$, the same arguments go
through, showing that $|D| \leqslant \ell^2$ or $|D| \leqslant
4\ell^2$.  One extra case is now possible, though: $y=0$.  This is the
case where $\p_\ell$ is \emph{ramified}, that is, $\p_\ell$ is not
principal but $\p_\ell^2 = (\ell)$.  In the ramified case, we must
have $\ell\mid D$.  In the rest of the proof, we will show that this
last case cannot in fact occur.  We state this as

%% Proof that $\ell\mid D$ is necessary.  Since $\p_\ell$ has
%% norm~$\ell$, there must exist a surjective ring homomorphism $\O
%% \to \F_\ell$. Such a map is determined by the image of $f\sqrt{d}$
%% or $f \frac{1+\sqrt{d}}{2}$ depending on the class of~$d$ mod~4.
%% Since~$\ell$ is odd, the ring homomorphism exists if and only if
%% $f^2d$ has a square root in~$\F_\ell$.  Since $\p_\ell$ is
%% ramified, the square root must be~0 in $\F_\ell$.  Thus $\ell\mid
%% f$ or $\ell\mid d$.  Hence $\ell\mid D$.]}

\begin{lemma} \label{noRam22}
With the hypotheses of the Theorem, if $\ell$ is odd, then $\ell \nmid
D$.
\end{lemma}

\begin{proof}
  %% Write $D = f^2 d_0$ as usual.  First, we show $\ell$ is relatively
  %% prime to the conductor~$f$.  In the ideal class of~$M$, choose a
  %% representative $M'$ so that both $M'$, and the superlattice $L'$
  %% corresponding to~$L$, are $\O$-ideals contained in~$\O$, and so
  %% that~$M'$ is relatively prime to~$f$.  Since $L'$ is an ideal that
  %% is a divisor of $M'$, $L'$ is also relatively prime to~$f$.  The
  %% quotient $\O$-ideal $L' \div M'$ is thus also relatively prime
  %% to~$f$ by [cite Cox or Cohn].  But this ideal has norm~$\ell$
  %% in~$\O$.  By [cite Cox again], $\ell$ is relatively prime to~$f$.
  %% Since $\p_\ell$ also has norm~$\ell$, it is also relatively prime
  %% to~$f$.
  %% Actually, what I'm about to say works with~$D$, so we don't need the
  %% separate argument for prime-to-$f$.
  Since $\p_\ell \not\sim (1)$, $L$ and $M$ are in different proper
  ideal classes.

  \cite[Prop.~7.3]{Bu} shows that when $\ell\mid D$, then among the
  $\ell+1$ index~$\ell$ sublattices of~$L$, there is precisely one
  that has the same coefficient ring.  This must be~$M$.  By the same
  argument, there is precisely one index~$\ell$ sublattice of~$M$ with
  the same discriminant, and this must be $\ell L$.

  In fact, however, there are two index~$\ell$ sublattices between $L$
  and $\ell L$.  This is because~$L$ and~$M$ have different ideal
  classes, hence different angles between their minimal vectors.
  $\tau \ne \sqrt{\ell}$ in this situation.  Of the two 2-and-2 forms
  $(L,M)$ and $(M,\ell L)$, one of them has $\tau < \sqrt{\ell}$ and
  the other $\tau > \sqrt{\ell}$.  The contradiction shows $\ell\mid
  D$ is impossible.
\end{proof}

When $\ell=2$, the analogue of Lemma~\ref{noRam22} still holds and is
proved by a direct computation.  This concludes the proof of
Theorem~\ref{thm:finManyDForEll}.
\end{proof}

\begin{remark}
  The bound $|D| \leqslant 4\ell^2$ is apparently not sharp.  When we
  list all the tempered perfect forms for $\ell < 100$, the largest
  ratio $\frac{|D|}{\ell^2}$ occurs at $\ell=47$, with $D = -6435$, so
  that $|D| \approx 2.91 \ell^2$.  We conjecture that $|D| \leqslant
  3\ell^2$ in general.
\end{remark}

\subsection{An algorithm for 2-and-2 forms}

As a complement to Algorithm~\ref{alg:eis}, we may combine
Corollary~\ref{oldconj5} and Theorem~\ref{thm:finManyDForEll} to give
an algorithm for finding all 2-and-2 forms whose index is a given
prime~$\ell$.  Let~$D$ run through all the negative discriminants
(integers $\equiv 0,1 \bmod 4$) from $-3$ through $-4\ell^2$.  As in
Lemma~\ref{noRam22}, omit~$D$ with $\ell\mid D$.  The prime~$\ell$
must split, so take only those~$D$ with $\left(\frac{D}{\ell}\right) =
+1$.  Find all reduced forms of discriminant~$D$; these are in
one-to-one correspondence with the proper ideal classes of the
order~$\O$ of discriminant~$D$.  Identify the well-rounded classes
$(a,b,a)$; if there are none, skip this~$D$.  Identify the
class~$\mathcal{P}$ that contains $\p_\ell$; it is given by the unique
reduced form that properly represents~$\ell$.  (The prime
$\overline{\p_\ell}$ is also in~$\mathcal{P}$, since~$\mathcal{P}$
must be an ambiguous class as in the proof of
Theorem~\ref{thm:finManyDForEll}.)  Whenever two well-rounded classes
$\mathcal{C}, \mathcal{C}'$ satisfy $\mathcal{C}\mathcal{P} =
\mathcal{C}'$, we get a 2-and-2 form, and all 2-and-2 forms of
index~$\ell$ arise in this way.

\subsubsection{Example: \textit{L} and \textit{M} with different shapes} \label{subsec:ex1155}

\begin{figure}
  \begin{center}
    \includegraphics[scale=0.19]{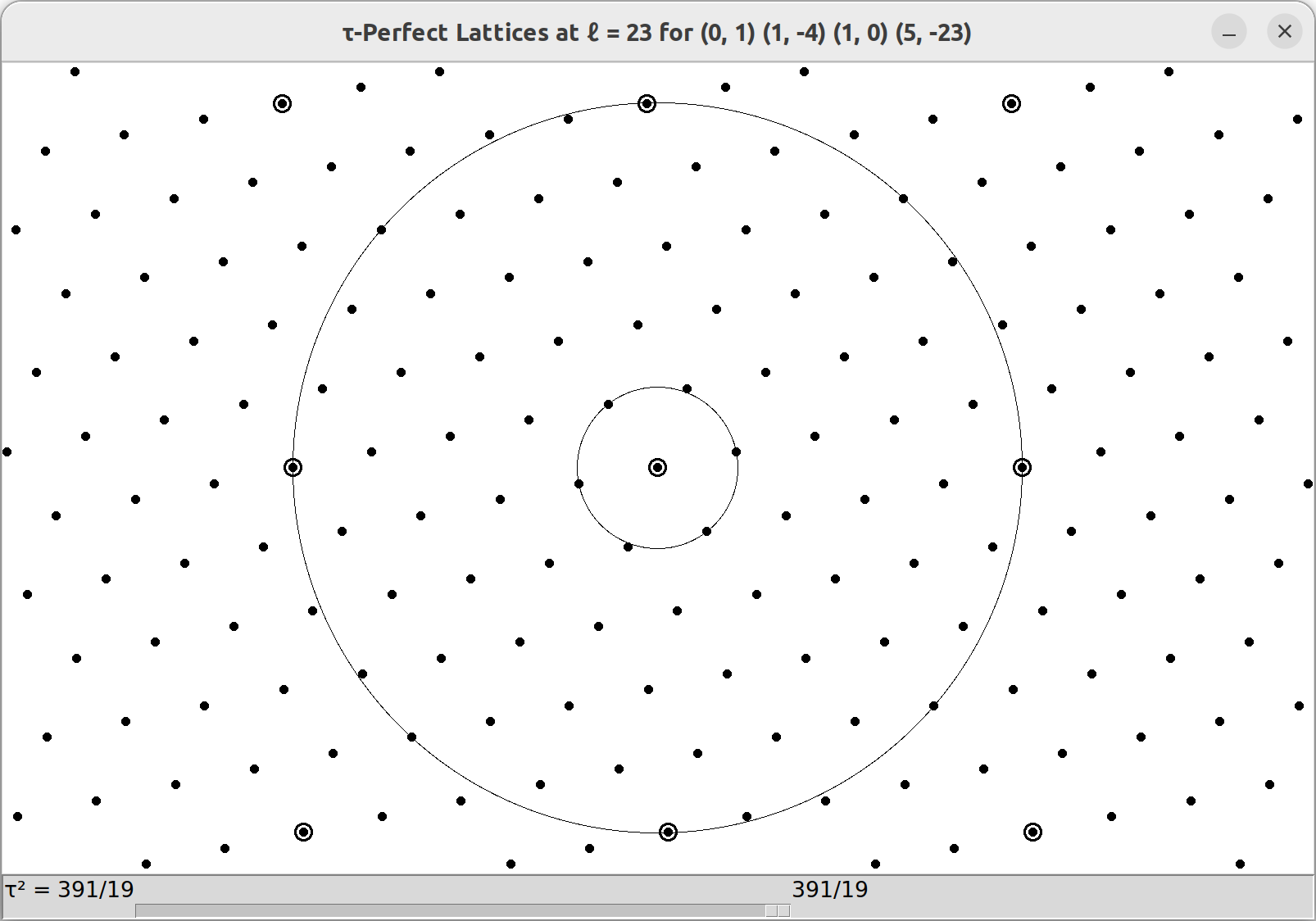}
  \end{center}
  \caption{A 2-and-2 form where $L$ and $M$ come from different ideal classes in a class group of order eight.}
  \label{fig:ex1155}
\end{figure}

At~$\ell = 23$, we have a 2-and-2 form $391x^2+169xy+19y^2$ on~$L_0 =
\Z^2$.  The discriminant $D = -1155 = -3\cdot 5 \cdot 7 \cdot 11$.
Let~$M_0$ be the $\Z$-span of $(1,0)$ and $(0,23)$.  The minimal
vectors in $L_0 - M_0$ are $(0,1)$ and $(1,-4)$, where the form takes
the value~19.  The minimal vectors in $M_0$ are $(1,0)$ and $(5,-23)$,
where the form takes the value $391 = 17\cdot 23$.  These are the
vectors on the inner and outer circles in Figure~\ref{fig:ex1155}.
The ratio of the circles' radii is $\tau = \sqrt{\frac{17\cdot
    23}{19}}$.

Let $\beta = \frac{1 + \sqrt{-1155}}{2}$.  The maximal order of $K =
\Q(\sqrt{-1155})$ is $\O_K = \Z[\beta]$.  The ramified primes are
$3,5,7,11$.  An odd prime~$p$ splits in $\O_K$ if and only if
$\left(\frac{-1155}{p}\right) = +1$.

The reduced forms of discriminant~$-1155$ appear in the first column
of Table~\ref{tab:1155}.  By~\eqref{bqfAmbig}, all eight forms are
ambiguous.  The class group is $\cong(\Z/2\Z)^3$.
Lemma~\ref{abaIsWR} gives a second proof that the forms
representing~17 and~19 are well rounded.  Only these two out of eight
classes are well rounded.

The ramified primes $\p_3$, $\p_5$, $\p_7$, and $\p_{11}$ each have
order two in the class group, and $\p_3 \p_5 \p_7 \p_{11} \sim (1)$ is
the only non-trivial relation among them.  Thus the eight reduced
forms correspond in some order to
$\p_3^{\varepsilon_3}\p_5^{\varepsilon_5}\p_7^{\varepsilon_7}$ with
exponents~$\varepsilon_p \in \{0,1\}$.  Calculation gives the middle
column of Table~\ref{tab:1155}.  The last column gives the smallest
non-inert rational prime~$p$ in the class.  (The~$331$ will be
explained in Section~\ref{subsubsec:ex1155genus}.)

\begin{table}
\begin{center}
  \begin{tabular}{|c|c|c|}
    \hline
    reduced form & generators & $p$ \\
    \hline
    \hline
    (1,-1,289) & $(1)$ & 331 \\
    \hline
    (3,-3,97) & $\p_3$ & 3 \\
    \hline
    (5,-5,59) & $\p_5$ & 5 \\
    \hline
    (7,-7,43) & $\p_7$ & 7 \\
    \hline
    (11,-11,29) & $\p_3\p_5\p_7$ & 11 \\
    \hline
    (17,-1,17) & $\p_5\p_7$ & 17 \\
    \hline
    (19,-17,19) & $\p_3\p_7$ & 19 \\
    \hline
    (15,-15,23) & $\p_3\p_5$ & 23 \\
    \hline
  \end{tabular}
\end{center}
\caption{The class group of $\Z[\sqrt{-1155}]$.}
\label{tab:1155}
\end{table}

The 2-and-2 forms of index~$\ell$ will arise from $D=-1155$ in two ways.
First, if~$\ell$ splits as a \emph{principal} ideal in~$\O_K$, then
$\p_\ell \cdot \p_5\p_7 \sim \p_5\p_7$ gives a 2-and-2 form where~$L$
and~$M$ both have the shape of the prime ideal over~17, and $\p_\ell
\cdot \p_3\p_7 \sim \p_3\p_7$ gives a 2-and-2 form where~$L$ and~$M$
both have the shape of the prime ideal over~19.

Second, note that $\p_3\p_5\cdot\p_5\p_7 \sim \p_3\p_7$.  Thus if
$\p_\ell \sim \p_3\p_5$, then multiplication by $\p_\ell$ will
interchange the two well-rounded ideal classes over~17 and~19.  This
is what is happening in the present example, where $\p_{23}\p_{19}
\sim \p_{17}$.

We emphasize that~$\p_{23}\p_{19} \ne \p_{17}$.  What matters is that
$\p_{23}\p_{19}$ has the same shape as~$\p_{17}$ up to rotation and
homothety.  The~$\p_{19}$ and~$\p_{17}$ are in different classes, but
both classes are well rounded.  The figure shows that the classes are
different.  The angle between the dots on the inner circle, a tiny bit
above $60^\circ$, is different from the angle between the~$\odot$
points on the outer circle, not far from~$90^\circ$.

In this second case, multiplying by $\p_\ell$ a second time would
perform the dual operation, $\p_\ell\p_{17} \sim \p_{19}$.  This
duality that exchanges two well-rounded but dissimilar classes was at
the heart of the proof of Lemma~\ref{noRam22}.

%%%%%%%%%%%%%%%%%%%%%%%%%%%%%%%%%%%%%%%%

\section{Two-and-two forms from a given order} \label{sec:22asFuncD}

%% Sources for the \emph{genus} of binary quadratic forms are
%% \cite[Ch.~14, 15]{Cohn} \cite[Chs.~3, 5, 6]{Cox}.

\subsection{Which indices appear for a given discriminant}

\begin{theorem} \label{thm:infManyEll}
  Let~$\O$ be the imaginary quadratic order of discriminant~$D = f^2
  d_0$ and conductor~$f$.  In the proper ideal class group of~$\O$,
  suppose $\mathcal{C}$ and $\mathcal{C}'$ are well-rounded ideal
  classes.  Then for all rational primes~$\ell$ for which
  \begin{equation} \label{inQuoClass}
    \mathrm{a\ proper\ prime\ ideal\ } \p_\ell
    \mathrm{\ of\ } \O \mathrm{\ is\ in\ the\ class\ }
    \mathcal{C}^{-1} \mathcal{C}',
  \end{equation}
  there exists a $2$-and-$2$ form for~$\ell$ in which the lattices
  $M\subset L$ are respectively in the classes $\mathcal{C}'$ and
  $\mathcal{C}$.  Infinitely many primes~$\ell$ satisfy
  condition~\eqref{inQuoClass}.

  In the special case where~$\O$ has only one ideal class per genus,
  there is a non-empty set of congruence classes modulo~$D$ such
  that~\eqref{inQuoClass} holds if and only if~$\ell$ is in one of
  those congruence classes (apart from a finite number of~$\ell$ that
  are exceptions).
\end{theorem}

\begin{proof}
  By Corollary~\ref{oldconj5}, if a prime $\p_\ell$ is in
  $\mathcal{C}^{-1} \mathcal{C}'$, then there exists a 2-and-2 form
  for~$\ell$.  By the Chebotaryov density theorem, the class
  $\mathcal{C}^{-1} \mathcal{C}'$ contains infinitely many prime
  ideals.

  When there is one class per genus, the statement about congruence
  classes follows from genus theory as developed by Lagrange and Gauss
  \cite[Thm.~2.26]{Cox}.  The finite number of exceptional~$\ell$ are
  those that divide~$D$.  We use only ideals of~$\O$ that are
  relatively prime to~$f$.  See \cite[(14.43)]{Cohn} when $f=1$, and
  see \cite[Exer.~18.16]{Cohn} or \cite[7.C]{Cox} when $f > 1$.
\end{proof}

\subsubsection{An example with one class per genus} \label{subsubsec:ex1155genus}

We continue the example of Section~\ref{subsec:ex1155}.  The genus
field \cite[(15.30)]{Cohn} over $K = \Q(\sqrt{-1155})$ is $L =
\Q(\sqrt{-3}, \sqrt{5}, \sqrt{-7}, \sqrt{-11})$.  The extension $L/\Q$
is Galois with group $(\Z/2\Z)^4$, and $\Q \subset K \subset
L$.  By class field theory, the class number of~$K$ will be a multiple
of~8, the degree of $L/K$.  Because the class number is exactly~8, the
genus field equals the Hilbert class field.  This implies
\cite[(15.33)]{Cohn} that~$p$ splits into a principal ideal of $\O_K$
if and only if
\[
\left(\frac{-3}{p}\right) = \left(\frac{5}{p}\right) =
\left(\frac{-7}{p}\right) = \left(\frac{-11}{p}\right) = 1.
\]
These four Legendre symbols together give the congruences mentioned in
Theorem~\ref{thm:infManyEll}.

For any discriminant $D<0$, the principal genus consists of the
squares in the class group.  Since the class group for $D=-1155$ has
exponent two, each genus contains exactly one class.  By Lagrange's
genus theory \cite[2.C]{Cox}, the $\frac12 \varphi(1155) = 240$ values
in $(\Z/1155\Z)^\times$ with $\left(\frac{D}{\dots}\right) = +1$ are
partitioned into eight cosets of 30 values each, each coset being the
set of values in $(\Z/1155\Z)^\times$ which the form can attain.

Our main result for $D=-1155$ is that Hecke operators for~$\ell$ will
have a 2-and-2 form with this~$D$ if and only if~$\ell$ is in the
values mod~1155 output by one of the two genera $(1,-1,289)$ and
$(15,-15,23)$, which correspond to the classes~$(1)$ and $\p_{23}$.
These values are respectively
\begin{gather*}
  1,4,16,64,169,214,256,289,331,361,379,394,421,466,499,526,529,631,
  \\ 676,694,709,751,841,856,949,961,991,1024,1054,1114; \\
  23,53,92,113,137,158,212,218,302,317,323,368,422,443,452,533,548,
  \\ 617,632,653,683,848,863,872,947,977,1037,1082,1103,1142.
\end{gather*}
The first few~$\ell$ in these two genera are 23, 53, 113, 137, 317 (in
the non-principal genus), 331, 379, 421 (principal genus), 443
(non-principal genus), etc.

\subsubsection{An example with two classes per genus}

Consider the discriminant $D = -55$.  The class group is $\Z/4\Z$.
The reduced forms are $(1,-1,14)$, $(2,-1,7)$, $(2,1,7)$, and
$(4,-3,4)$.  The squares in the class group are $(1,-1,14)$ and
$(4,-3,4)$, so these two constitute the principal genus.  The values
mod~55 taken by the forms in the principal and non-principal genera
are respectively
\[
\underbrace{1,4,9,14,16,26,31,34,36,49}_\mathfrak{g}
\quad\mathrm{and}\quad 2,7,8,13,17,18,28,32,43,52.
\]
Since there is only one well-rounded reduced form, $(4,-3,4)$,
the~$\ell$ we need are those in the principal class, because they
carry the well-rounded class to itself.

The difficulty, compared to the last example, is that congruence
conditions on~$\ell$ are not enough to pick out the principal class.
If $\ell\in\mathfrak{g} \bmod{55}$, then we have no way of knowing
from this congruence information alone whether $\p_\ell$ is in the
principal class or in the $(4,-3,4)$ class.

To distinguish our~$\ell$, we may use the Hilbert class field.  The
Hilbert class field of $K = \Q(\sqrt{-55})$ is a $\Z/4\Z$ extension
of~$K$ and is a quadratic extension of the genus field $\Q(\sqrt{5},
\sqrt{-11})$.  Calculations show that the Hilbert class field is the
splitting field of $x^4 - x^3 + 2x - 1$ over~$\Q$ and that the four
roots of this polynomial are
\[
\varepsilon_1 \frac14 \sqrt{5} + \varepsilon_2 \frac12
\sqrt{-\varepsilon_1 \frac32 \sqrt{5} - \frac12} + \frac14
\]
where $\varepsilon_1 = \pm1$ and $\varepsilon_2 = \pm1$ independently.
See \cite[Exercise~6.20]{Cox}.

To determine in practice whether a given $\ell\in\mathfrak{g}
\bmod{55}$ is in the principal class, we need not use the Hilbert
class field directly.  We can simply check whether $x^2-xy+14y^2$ or
$4x^2-3xy+4y^2$ represents~$\ell$; it will be one or the other.  For
instance,
\begin{align*}
  31 &= 4(1)^2 - 3(1)(3) + 4(3)^2 \qquad\mathrm{(nonprincipal\ class)} \\
  59 &= (1)^2 - (1)(-2) + 14(-2)^2 \qquad\mathrm{(principal\ class)} \\
  71 &= (3)^2 - (3)(-2) + 14(-2)^2 \qquad\mathrm{(principal\ class)} \\
  89 &= 4(1)^2 - 3(1)(5) + 4(5)^2 \qquad\mathrm{(nonprincipal\ class)}.
\end{align*}
The~$\ell$ where we see 2-and-2 forms for $D=-55$ are $\ell=59$, $71$,
etc.

\subsection{Doubly well-rounded forms have round discriminants} \label{subsec:ruleOfThree}

When an order~$\O$ has a well-rounded ideal class, we observe that its
discriminant~$D$ is a product of small primes.  When~$\O$ has more
than one well-rounded class, $D$ is even more strikingly a product of
small primes.  The first examples of~$D$ with two different
well-rounded classes are $-1155 = -3\cdot 5 \cdot 7 \cdot 11$ and
$-1120 = (4)(-2^3\cdot 5 \cdot 7)$.

Let $(a,b,c)$ be a reduced form of discriminant $D<0$.  By
Lemma~\ref{abaIsWR}, the form is well rounded if and only if $a=c$,
implying $D = b^2 - 4a^2$.  Thus we would like to know for which
discriminants $D<0$ we can solve
\begin{equation} \label{qu2}
D = b^2 - 4a^2 \quad\mathrm{with\ } |b|\leqslant a \mathrm{\ and\ } (a,b)=1,
\end{equation}
where the last two conditions are being reduced and primitive.  A
solution to~\eqref{qu2} guarantees by Theorem~\ref{thm:infManyEll}
that there are 2-and-2 forms coming from discriminant~$D$ for
infinitely many~$\ell$, those in the identity class.  Having two or
more solutions guarantees there are 2-and-2 forms $(L,M)$ for
infinitely many~$\ell$ in which~$L$ and~$M$ come from different
classes.

We first suppose~$D$ is odd.

\begin{proposition}
If $D\equiv 1 \pmod 4$, there is a solution to~\eqref{qu2} if and only
if~$D$ factors as $-FG = D$ with $F,G$ relatively prime and satisfying
\[
\frac13 \leqslant \frac{F}{G} \leqslant 3.
\]
\end{proposition}

\begin{proof}
  We have a solution if and only if $-(2a-b)(2a+b) = D$ with
  $|b|\leqslant a$ and $\gcd(a,b)=1$.  Write $-FG = D$ with $F = 2a -
  b$ and $G = 2a+b$.  Clearly $F,G\geqslant 1$, $\frac{F+G}{4} = a$,
  and $\frac{G-F}{2} = b$.  From $|b|\leqslant a$ we obtain
  $-\frac{F+G}{4} \leqslant \frac{G-F}{2} \leqslant \frac{F+G}{4}$.
  By algebra, this gives $\frac13 \leqslant \frac{F}{G} \leqslant 3$.

Secondly, $\gcd(a,b) = \gcd\left(\frac{F+G}{4}, \frac{G-F}{2}\right) =
\gcd\left(\frac{F+G}{4}, G\right)$ by adding twice the first argument
to the second.  Both $F,G$ are odd because their product is odd.
Because~$G$ is odd, the gcd is unchanged if we multiply the first
argument by a power of~2.  Hence the gcd is $\gcd(F+G, G) =
\gcd(F,G)$.
\end{proof}

\begin{examples}
  $D=-3$ has a well-rounded lattice (the hexagonal lattice) because
  $\frac31 = 3$.  This shows the bound in the proposition is sharp.

  $D=-55$ has a well-rounded lattice because $\frac{11}{5} = 2.2 < 3$.
  
  The first case with two solutions is $D=-1155 =
  -3\cdot5\cdot7\cdot11$.  Here $\frac{35}{33}$ is just a little
  above~1, and $\frac{55}{21}$ is still below~3.
\end{examples}

Here is a similar result for even~$D$.

\begin{proposition}
If $D\equiv 0 \pmod 4$, there is a solution to~\eqref{qu2} if and only
if~$D$ factors as $-FG = D$ with $F,G$ satisfying
\[
\frac13 \leqslant \frac{F}{G} \leqslant 3
\quad\mathrm{and}\quad
\gcd\left(\tfrac{F+G}{2}, 2G\right) = 1.
\]
\end{proposition}

\begin{proof}
Here $b$ must be even.  Write $b = 2b_1$, so $b_1^2 - a^2 = D/4$.  We
have a solution to~\eqref{qu2} if and only if $-(a-b_1)(a+b_1) = D/4$
with $|b_1| \leqslant a/2$ and $\gcd(a,2b_1) = 1$.  Write $-FG = D/4$
with $F = a-b_1 \geqslant 1$ and $G = a+b_1 \geqslant 1$.  Solving,
$\frac{F+G}{2} = a$ and $\frac{G-F}{2} = b_1$.  From $|b_1| \leqslant
a/2$ we deduce the same formulas $-\frac{F+G}{4} \leqslant
\frac{G-F}{2} \leqslant \frac{F+G}{4} \Rightarrow \frac13 \leqslant
\frac{F}{G} \leqslant 3$ as for odd~$D$.
\end{proof}

%% [The condition $\gcd(a,b)=1$ seems to break into a lot of cases.
%%   $D/4$ can have any congruence class mod~4.  Can $\gcd(a,b)=1$ be
%%   expressed cleanly as one or a few conditions on $\gcd(F,G)$,
%%   $\gcd(F/2, G/2)=1$, etc.?  Here is a start.]  $F,G$ have the same
%% parity mod~2, as $a\in\Z$ implies.  $\gcd\left(\frac{F+G}{2},
%% G-F\right) = \gcd\left(\frac{F+G}{2}, 2G\right)$ by adding twice the
%% first argument to the second.  A necessary condition for
%% $\gcd\left(\frac{F+G}{2}, 2G\right) = 1$ is $F+G\equiv 2 \pmod 4$.

\begin{examples}
  $D=-4$ with $\frac11 < 3$ and $\gcd(1,4)=1$.

  $D=-32$ with $\frac42 < 3$ and $\gcd(3,4)=1$.

  The first case with two solutions is $D=-1120 = (4)(-280)$, where
  $-280 = -2^3 \cdot 5 \cdot 7$.  Here $\frac{20}{14} = 1.429 < 3$
  with $\gcd(17, 28) = 1$, and $\frac{28}{10} = 2.8 < 3$ with
  $\gcd(19, 20) = 1$.
\end{examples}

%%%%%%%%%%%%%%%%%%%%%%%%%%%%%%%%%%%%%%%%

\bibliographystyle{plain}
\bibliography{tauperf}

\begin{thebibliography}{10}

\bibitem{Ash84}
Avner Ash.
\newblock Small-dimensional classifying spaces for arithmetic subgroups of general linear groups.
\newblock {\em Duke Math.\ J.}, 51(2):459--468, June 1984.

\bibitem{AGMY}
Avner Ash, Paul Gunnells, Mark McConnell, and Dan Yasaki.
\newblock On the growth of torsion in the cohomology of arithmetic groups.
\newblock {\em J.\ Inst.\ Math.\ Jussieu}, 19(2):537--569, 2020.

\bibitem{AR}
Avner Ash and Lee Rudolph.
\newblock The modular symbol and continued fractions in higher dimensions.
\newblock {\em Invent.\ Math.}, 55:241--250, 1979.

\bibitem{BoSh}
Z.~I. Borevich and I.~R. Shafarevich.
\newblock {\em Number Theory}.
\newblock Academic Press, 1966.

\bibitem{Bu}
Duncan~A. Buell.
\newblock {\em Binary Quadratic Forms: Classical Theory and Modern Computations}.
\newblock Springer-Verlag, 1989.

\bibitem{Cohn}
Harvey Cohn.
\newblock {\em A Classical Invitation to Algebraic Numbers and Class Fields}.
\newblock Springer-Verlag, 1978.

\bibitem{CS}
J.~H. Conway and N.~J.~A. Sloane.
\newblock {\em Sphere Packings, Lattices and Groups}, volume 290 of {\em Grundlehren der math.\ Wiss.}
\newblock Springer-Verlag, third edition, 1999.

\bibitem{Cox}
David~A. Cox.
\newblock {\em Primes of the Form $x^2 + ny^2$}.
\newblock Wiley, second edition, 2013.

\bibitem{GM21}
Dylan Galt and Mark McConnell.
\newblock Computing {H}ecke operators for arithmetic subgroups of {S}p${}_4$.
\newblock {\em Advanced Studies: Euro-Tbilisi Math.\ J.}, pages 67--80, Sept. 2021.

\bibitem{Gun}
Paul Gunnells.
\newblock Computing hecke eigenvalues below the cohomological dimension.
\newblock {\em J.\ Experimental Math.}, 9(3):351--367, 2000.

\bibitem{GunMY}
Paul Gunnells, Mark McConnell, and Dan Yasaki.
\newblock Cohomology of {GL}(3) over the eisenstein integers.
\newblock {\em J.\ Experimental Math.}, 30(4):499--512, 2021.

\bibitem{MM89}
Robert MacPherson and Mark McConnell.
\newblock Classical projective geometry and modular varieties.
\newblock In {\em Algebraic Analysis, Geometry, and Number Theory}, pages 237--290. Johns Hopkins U.\ Press, 1989.

\bibitem{Manin}
Y.~I. Manin.
\newblock Parabolic points and zeta-functions of modular curves.
\newblock {\em Math.\ USSR Izvestija}, 6(1):19--64, 1972.

\bibitem{MM20}
Mark McConnell and Robert MacPherson.
\newblock Computing {H}ecke operators for arithmetic subgroups of general linear groups.
\newblock arXiv math.NT 2010.06036, 2020.

\bibitem{NZM}
Ivan Niven, Herbert~S. Zuckerman, and Hugh~L. Montgomery.
\newblock {\em An Introduction to the Theory of Numbers}.
\newblock Wiley, fifth edition, 1991.

\bibitem{Stein}
William Stein.
\newblock {\em Modular Forms: A Computational Approach}, volume~79 of {\em Grad.\ Studies in Math.}
\newblock Amer.\ Math.\ Soc., 2007.

\bibitem{Vor}
Georges Voronoi.
\newblock Nouvelles applications des param\`{e}tres continus \`{a} la th\'{e}orie des formes quadratiques: 1~{S}ur quelques propri\'{e}t\'{e}s des formes quadratiques positives parfaites.
\newblock {\em J.\ Reine\ Angew.\ Math.}, 133:97--178, 1908.

\end{thebibliography}

\end{document}